\documentclass{amsart}
\usepackage{geometry}
\usepackage{graphicx}
\usepackage{amssymb}
\usepackage{natbib}
 
\DeclareMathOperator*{\argmax}{arg\,max} 
\def\Halmos{}%\fbox{\phantom{\rule{.7ex}{.7ex}}}
\newtheorem{theorem}{Theorem}
\newtheorem{corollary}{Corollary}
\newtheorem{lemma}{Lemma}
\newtheorem{proposition}{Proposition}
\newtheorem{remark}{Remark}

\newcommand{\indep}{\rotatebox[origin=c]{90}{$\models$}}
\newcommand{\mC}{{\mathcal C}}

\newcommand{\mJ}{{\mathcal J}}
\newcommand{\mK}{{\mathcal K}}
\newcommand{\mL}{{\mathcal L}}

\newcommand{\mR}{{\mathcal R}}
\newcommand{\mS}{{\mathcal S}}

\newcommand{\bN}{{\mathbb N}}
\newcommand{\bR}{{\mathbb R}}
\newcommand{\bZ}{{\mathbb Z}}
\newcommand{\bP}{{\mathbb P}}
\newcommand{\bE}{{\mathbb E}}

\newcommand{\vA}{{\mathbf A}}

\newcommand{\vM}{{\mathbf M}}

\newcommand{\vQ}{{\mathbf Q}}

\newcommand{\ve}{{\mathbf e}}

\newcommand{\vm}{{\mathbf m}}

\newcommand{\vq}{{\mathbf q}}

\newcommand{\coS}{<\!{\mathcal S}\! >}

\newcommand{\out}{\omega}
\newcommand{\invA}{A^{-1}}
\newcommand{\jInJ}{{j \in \mJ}}
\newcommand{\rInR}{{r \in \mR}}

\newcommand{\kInj}{{k \in j}}

\newcommand{\vecA}{{\bf A}}
\newcommand{\vecD}{{\bf D}}
\newcommand{\vecQ}{{\bf Q}}
\newcommand{\vecX}{{\bf X}}
\newcommand{\vecpi}{{\boldsymbol \pi}}

\newcommand{\vecalpha}{{\boldsymbol \alpha}}

\newcommand{\vecGamma}{{\boldsymbol \Gamma}}

\newcommand{\veca}{{\boldsymbol a}}
\newcommand{\vecp}{{\bf p}}
\newcommand{\vecq}{{\bf q}}

\newcommand{\sigmamax}{\sigma_{\rm{max}}}
\newcommand{\amin}{a_{\rm{min}}}
\newcommand{\vecsigma}{{\boldsymbol \sigma}}

\begin{document}

% Title
\title[Proportional Switching]{Proportional switching in FIFO networks}

% Authors
% first
\author{Maury Bramson}
\address{School of Mathematics, University of Minnesota, 206 Church St. SE, MN 55455, USA.}
\email{bramson@math.umn.edu}
\thanks{The research of the first author was partially supported by NSF grants DMS-1105668 and DMS-1203201.}
% second
\author{Bernardo D'Auria}
\address{Departmento de Estad\'istica, Universidad Carlos III de Madrid, Calle Madrid 126, 28903 Getafe, Madrid, Spain.}\email{bernardo.dauria@uc3m.es}
\thanks{The research of the second author was partially supported by the Spanish Ministry of Economy and Competitiveness Grants MTM2013-42104-P via FEDER funds;
he  thanks the ICMAT (Madrid, Spain) Research Institute that kindly hosted him while developing this project.}
% third
\author{Neil Walton}
\address{Korteweg-de Vries Institute for Mathematics, University of Amsterdam, Science Park 904, 1098 XH Amsterdam, NL.}\email{n.s.walton@uva.nl}
\thanks{The research of the first author was partially supported by NSF grants DMS-1105668 and DMS-1203201.
The research of the second author was partially supported by the Spanish Ministry of Economy and Competitiveness Grants MTM2013-42104-P via FEDER funds;
he  thanks the ICMAT (Madrid, Spain) Research Institute that kindly hosted him while developing this project.
The research of the third author was funded by the VENI research programme, which is financed by the
Netherlands Organisation for Scientific Research (NWO)}


\date{}

%\dedicatory{}

\begin{abstract}
We consider a family of discrete time multihop switched queueing networks where each packet moves along a fixed route. In this setting, BackPressure is the canonical choice of scheduling policy; this policy has the virtues of possessing a maximal stability region and not requiring explicit knowledge of traffic arrival rates.  BackPressure has certain structural weaknesses because implementation requires information about each route, and queueing delays can grow super-linearly with route length. For large networks, where 
packets over
many routes are processed by a queue, or where packets over a route are processed by many queues, these limitations can be prohibitive. 
%A well-known alternative would be the MaxWeight policy; however, MaxWeight does {not} have a maximal stability region in this setting.
%As is common practice in numerous Internet applications, we aggregate packets into first-in, first-out queues. 

In this article, we introduce a scheduling policy for FIFO networks, the Proportional Scheduler, which is based on the proportional fairness criterion. We show that, like BackPressure, the Proportional Scheduler has a maximal stability region and does not require explicit knowledge of traffic arrival rates.
The Proportional Scheduler has the advantage that information about the network's route structure  is not required for scheduling, 
which substantially improves the policy's performance for large networks.
For instance, packets can be routed with only next-hop information and new nodes can be added to the network with only knowledge of the scheduling constraints. 
%
%In addition, the Proportional Scheduler shares some interesting relationships with the classical quasi-reversible queueing networks; after briefly discussing them we show why these would imply, for the former, a desirable \emph{product form resource pooling} property as well as nice delay scaling effects.
\end{abstract}

% MAKETITLE
\maketitle

% SECTIONS
\section{Introduction}
% EXMAPLES TO HELP MAURY
%\red{added text} \cancel{canceled text}
% DESCRIBE A SWITCHED NETWORK 
We consider, in this paper, a family of discrete time multihop switched queueing networks where packets move along fixed routes. Switched networks were first introduced in \cite{TaEp92}; this
terminology was first employed in \cite{ShWi12} for a discrete time queueing network whose queues are served simultaneously subject to certain global scheduling constraints and with packets moving along fixed routes. 
Here, we consider a variant of this model that allows more than one route through
each queue; the scheduling constraints and queueing discipline employed here do not depend on the routes of the packets.

 Applications of switched networks include 
 wireless ad-hoc networks,
 Internet routers, call centers with cross trained staff,
 data centers, and
 urban road traffic scheduling. 
For such applications, and in general,
%
%
%For a switched policy, an important first order property is to have a maximal stability region.
%
maximum stability is a highly desirable feature for the scheduling policy.
Roughly stated, for switched networks, a scheduling policy is maximally stable if, for every arrival rate for which there exists a stable policy, the policy stabilizes the network.
When the vector of arrival rates is known, one can specify a stable policy
by choosing a random schedule whose average service rate at each queue dominates the corresponding arrival rate.
However, in practice, explicit knowledge of arrival rates is often not available, particularly when rates may vary. Since the seminal work of \cite{TaEp92},  BackPressure has been the canonical maximally stable scheduling policy 
for multihop switched networks. An important feature of this policy is that, in addition to it
being maximally stable, only information on the local state is required for a scheduling decision.

When packets are routed through different network components, policies such as BackPressure require detailed knowledge of this routing.  Such policies are similar in nature to the classical single class Jackson queueing networks if one identifies each route with a queue, so that each queue now contains only a single customer type. However, in many applications, this single class interpretation may not be practical  because the information required can increase rapidly in the number of queues of the original network.  It is also not common practice --
for instance, an Internet router will maintain a first-in, first-out (FIFO) queue for each outgoing link rather than a queue for each route, because the number of links to and from a router is orders of magnitude smaller than the number of routes it  processes. 
%In the analysis of switched networks, it is often assumed that queues are single class. Here each queue only contains one type of customer.
%
As a result, in many practical situations, the model simplifies if one reinterprets 
queues as being multiclass, i.e., permitting different packet classes and hence different routing of packets passing through the queue. 

In our setting, queues will be multiclass but will serve packets according to a scheduling policy allocating service among queues that depends only on the number of packets at each queue; 
the discipline at each queue will be FIFO and packets will all be of unit size. 
In this setting, we will show that this switching policy, the Proportional Scheduler, is maximally stable.
%
%Maximal stability has been shown when packets leave the network immediately after first being served, so called single hop networks.
%

There are well-known examples of disciplines that are not maximally stable 
for multiclass queues, for instance, \cite{LuKu91}, \cite{RySt92} and \cite{Br94}.  (These examples are in both discrete and continuous time.)  The examples in \cite{Br94} are for the FIFO discipline, with jobs having unequal mean service requirements.  In switch networks, all jobs have equal service requirements and so do not fall within this framework.

Our policy, the \emph{Proportional Scheduler}, can be described roughly as follows: for a set of FIFO queues $\mJ$, a vector of queue lengths $(Q_j: j\in\mJ)\in \bZ_+^{|\mJ|}$, and a convex set of schedules $ <\!\mS \!> \subset \bR^{|\mJ|}$, the Proportional Scheduler serves packets according to a vector of expected rates,  $\vecsigma \in \bR^{|\mJ|}$, that solves the proportional fair optimization problem
\begin{align}\label{IntroPF}
\text{maximize}\qquad \sum_{j\in\mJ} Q_j \log \sigma_j \qquad
\text{over}\qquad  \vecsigma \in <\! \mS \! >.
\end{align}
Packets within a queue are to be served according to a FIFO queueing discipline; once served, a packet goes to the next queue on its route.

We  prove that maximal stability holds for this policy by employing a fluid model analysis of our system. This combines the approaches of \cite{Ma07}, for bandwidth networks, and \cite{Br96a}, for FIFO Kelly networks; the latter networks can be viewed as the special case of switched networks where the allocation of
service to different queues is fixed irrespective of the state of the system.

We also compare our policy to the BackPressure policy. Unlike BackPressure, the Proportional Scheduler does not require knowledge of the route used by each packet. Thus, routing structure and scheduling decisions are distinct
for the Proportional Scheduler policy. 
%From a practical perspective, this policy does not require information on the routes of packets in order to make a scheduling decision.
As we later discuss in more detail, there are various benefits of the Proportional Scheduler's functional structure: 
a) rather than distinguishing packets according to their routes (either in the memory of the scheduling algorithm or by physically maintaining different queues), packets can be served at their queue in a first-in, first-out order, 
b) a packet can be routed knowing only its next hop, rather than knowing its entire route,
c) when adding new network components, one only needs to know the component's scheduling constraints and not the entire route or class structure of the network to implement the policy, and
d) for BackPressure, messages must be sent between queues, which the Proportional Scheduler does not require. 

On the other hand, the Proportional Scheduler policy is less general than BackPressure in that the latter model allows for adaptive routing decisions whereas routing for the Proportional Scheduler is fixed.

Another well-known switched network policy is the MaxWeight policy; it is defined by the optimization whose objective function employs the linear factor $\sigma_j$, in place of $\log \sigma_j$, in \eqref{IntroPF}, i.e.,
\begin{align}\label{IntroMW}
\text{maximize}\qquad \sum_{j\in\mJ} Q_j  \sigma_j \qquad
\text{over}\qquad  \vecsigma \in <\! \mS \! >.
\end{align}
  As such, it also satisfies the properties a) -- d) above.  MaxWeight is maximally stable for single hop networks, since it coincides there with the BackPressure policy.  

MaxWeight is not
always maximally stable for multihop networks, even when each queue possesses
only one class.  This is shown in \cite{andrews2003achieving} for
fluid models for a weighted version of MaxWeight, after modifying the model slightly by varying the arrival rate of mass; supporting simulations are then given for a queueing network under the standard MaxWeight policy.  
We briefly consider the maximal stability of MaxWeight in Section \ref{Comparison}, but leave a detailed investigation of 
the policy to future work.

The same question about maximal stability exists for policies with objective functions employing factors other than 
$\sigma_j$ or $\log \sigma_j$.  As the counterexample in \cite{andrews2003achieving} shows, care needs to be
exercised in the choice of the objective function in \eqref{IntroPF}.

\subsection{Relevant Literature}
The main results of this paper combine a number of results, methods, and models from the theory of queueing networks over the last three decades. We review this relevant literature  in a roughly chronological order.

% CLASSICAL PRODUCT FORM NETWORKS
\emph{A. Classical Queueing Networks:} The development of Jackson networks is one of the earliest substantial developments in the theory of stochastic networks.
In Jackson networks, a single class of customers is routed probabilistically between queues in the network.
Influenced by the Input Theorem of \cite{Bu56}, \cite{Ja63} found that the stationary distribution of these networks can be written in product form, meaning that its stationary distribution is the product of simple terms. 
%For Jackson networks, a single class of customer is routed probabilistically between queues in a network.
%
%As with BackPressure policies, these networks are single class and so sophistication of the networks routing structure can only be increased by increasing the number of queues within the network.
%
%However, numerous practical situations require a queue to be multitype, e.g., able to
%process packets passing through it that traverse different routes, which requires generalization of Jackson networks.
%

%
\cite{BCMP75} and \cite{Ke75,Ke79} 
broadened this family of queueing networks by permitting queues to be multiclass, and hence allowing more than one route through each queue. 
Using quasi-reversibility, they showed that, for certain service disciplines including FIFO, the stationary distribution must be of product form.  (Quasi-reversibility will 
also be
employed in the context of proportional fairiness in Part D below.)

%This analyse provided the first tractable analysis of networks of FIFO queueing networks. Their results exploited the 
%%
%%\NW{\sout{reversibility and}} 
%%
%quasi-reversiblity of certain queueing disciplines, \NW{which include first-in first-out queueing,} \cite{Mu72,Ke79}. 
%\NW{(Further reference to quasi-reversible queueing networks will be made when we review proportional fairness in Part D of this section.)}

%\NW{{DELETE THE FOLLOWING BLUE TEXT: Reversible queueing systems were extended and analysed by  \cite{Wh85}, and reversibility still remains one of the most tractable ways of analysing connected queueing components.  The quasi-reversibility property for Whittle networks motivates our proof of 
%maximum stability for proportional switching in multiclass networks. Relationships between our proportional fairness and quasi-reversible Whittle networks are described in more detail in \cite{Ma07}. Thus there is a close relationship between Proportional Scheduling and quasi-reversible queueing networks. For single-hop FIFO networks, product form heavy traffic results are derived by \cite{Ye12}  and, for a closely related store-forward allocation, some discussions of quasi-reversibility and product-form in the context of continuous time FIFO networks can be found in \cite{walton2014store}.}}

% BackPressure
\emph{B. Switched Networks, BackPressure and their Applications:} Switched queueing networks and, more specifically, the BackPressure scheduling policies, were first introduced by \cite{TaEp92} as a model of wireless communication. As mentioned above, a switched queueing network is a discrete time queueing network with constraints on which queues can be served simultaneously.
The BackPressure policies have proven popular because they maximize a network's stability region while not requiring explicit estimation of traffic arrival rates;
for a comprehensive review of the BackPressure policies, see \cite{GNT06}.

The BackPressure policies have been generalized and specialized in numerous directions. 
In contrast to our work, the defining feature of a BackPressure policy is that the policy minimizes the drift of a Lyapunov function subject to the scheduling constraints of the network. For single hop networks -- where packets are served only once before departing the network -- the BackPressure policy is often referred to as the MaxWeight policy.
As a model of Internet protocol routers, \cite{MMAW99} applied this paradigm to the example of input-queued switches.
 \cite{AKRSVW04} consider power functions when defining their MaxWeight Lyapunov function and thus generalized the set of MaxWeight/BackPressure policies.
Additional extensions are considered by \cite{Me09} and \cite{ESP05} and further generalizations to cone polices are considered by \cite{ArBa03}.
In different stochastic senses, BackPressure and MaxWeight can be shown to optimize certain workload functions:
 in heavy traffic, see \cite{St04}; in large deviations, see \cite{Ve07}; in overload, see \cite{ShWi11}. Further, heavy tailed arrivals are analyzed by \cite{JMMT11};
delay in the presence of heavy tailed traffic are considered in \cite{MMT14}.

 %
  % Applications of Switch networks
 
There are numerous application areas associated with switched networks. For these areas, BackPressure is often the canonical choice. Applications  include
 wireless ad-hoc networks, in \cite{TaEp92},
 Internet routers, in \cite{MMAW99},
 call centers with cross trained staff, in \cite{MaSt04},
 data centers, in \cite{ShWi11},
 urban road traffic scheduling, in \cite{Va13}, and
 stochastic processing networks, which would include numerous manufacturing and general processing settings, in \cite{DaLi05}.

% Instability, Fluid Stability and Stability
\emph{C. Instability, Fluid Stability and Stability of Queueing Networks:}
Most classical queueing networks are positive recurrent when the network is subcritical.  By positive recurrent, we mean that there exists a state (which, in our setting, will be the state with no packets) that is positive recurrent and that this state is visited with probability $1$ starting from any other state; this is the
reduction to the countable state space setting of positive Harris recurrence, which is the standard definition in the
general state space setting.  By subcritical, we mean that each network resource experiences a load that is strictly less than the resource's capacity. It had been thought that subcritical networks were always positive recurrent under a work conserving policy. However, a series of examples constructed in the mid-nineties showed that this is not the case.   For instance, see \cite{LuKu91}, \cite{RySt92} and \cite{Br94}.

This led to new approaches for determining the stability region for queueing networks. In particular, \cite{RySt92} and \cite{Da95} developed a fluid model approach where the stability of a queueing network can be determined from that of an associated fluid model. This theory is surveyed in by \cite{Br08}.
A fluid analysis of multiclass FIFO queueing networks was first given by \cite{Br96a}. Our fluid analysis uses a similar approach.

Recent work of \cite{DiSh13} further considers the issue of finding natural scheduling policies that lead to stability whenever each server is nominally underloaded. Similar to adversarial queueing frameworks, e.g., \cite{BKRSW01}, stability is achieved by prioritizing queue service according to a least-routed-first-priority discipline.  
This differs from the approach taken in this paper where the service discipline does not use the routing structure of the network  to achieve maximum stability.
Further recent work of \cite{JJS13} considers stability results for switched networks where packets are queued per-link rather than per-route, as is typically applied for BackPressure. Here stability is achieved by running an appropriate {queueing system} in the memory of the algorithm. By estimating queue sizes and loads in this way, stability is achievable for switched networks as long as routes do not form a loop. Once again this differs from the approach of this paper, where routes are general and only current queue size information is required to execute the scheduling policy. Finally, \cite{ying2011cluster} consider a modification of the BackPressure policy to lessen the effect of per-route queueing. Queues are grouped into predetermined clusters, with the standard BackPressure policy applied within each cluster. A judicious choice of clusters will reduce the amount of required memory. Clusters are centrally determined and, as with BackPressure, per-link queueing is not employed and information on the queue size must be continuously exchanged along queues between the source and destination.

% Bandwidth networks
\emph{D. Massouli\'e Networks, Proportional Fairness and Quasi-Reversibility:}
A further class of  Internet models was introduced in \cite{MaRo99} and
\cite{massoulie2000bandwidth}
in a processor-sharing framework where resources are shared subject to constraints on these resources. 
%
%\NW{Rather than switching between different discrete schedules, \cite{MaRo99} assume a processor-sharing framework where resources are shared subject to constraints on these resources.} 
%
Similar to MaxWeight and BackPressure, these policies are often defined by an optimization that is maximally stable.  (Unlike BackPressure, the construction does not minimize the instantaneous drift of a Lyapunov function.) Stability proofs for these systems can be found in \cite{BoMa01}, \cite{Ye05}, \cite{GW09} and \cite{PAFA12}. Due to their proliferation to different areas, these models have taken various names, such as bandwidth sharing networks, stochastic flow level models, and resource sharing networks; here, we refer to these networks as \emph{Massouli\'e networks}. See \cite{HMSY14} for a recent discussion of the benefits and varied applications of this resource sharing paradigm.

In this paper, we allocate resources according to the proportional fair optimization, 
which was first introduced by \cite{Ke97}. Proportional fairness has been 
used in the allocation of bandwidth in modern 3G telephone networks, see, e.g., \cite{VTL02} and \cite{KuWh04}.  The stability of proportional fairness in Massouli\'e networks was first shown by \cite{DLK99}; an important generalization of this stability analysis is given in \cite{Ma07}. Further progress on the stability and large deviations behavior of proportional fairness can be found in \cite{JoLo14}. Heavy traffic analysis of proportional fair policies can be found in \cite{KKLW07i}, \cite{Ye12}, and \cite{VZZ14}. \cite{St04} has investigated resource pooling for MaxWeight policies. \cite{KKLW07ii},  \cite{KKLW07i}, \cite{KMW09} and \cite{Ye12} discussed product form resource pooling properties associated with the proportional fairness in heavy traffic and large deviations regimes.

Influenced by \cite{Wh85}, the quasi-reversibility and insensitivity property in Massouli\'e networks was studied by \cite{BoPr02,BoPr03,BoPr04} and \cite{Za07}. 
For connections between proportional fairness and the queueing networks of \cite{Ke75} and \cite{BCMP75}, see, e.g., \cite{Sc79,Ke89,MaRo99,Wa09}, and \cite{MOR13}.
A short, general description of the relationship between proportional fairness, maximum stability and quasi-reversibility can also be found in \cite{Wa11B}.
 As we discuss later, the Lyapunov function in \cite{Ma07} is relevant to our analysis.

\emph{E. Resource Sharing in Switched Networks:}
The sharing of network resources is a key property of the proportional fair optimization. Other resource sharing policies exist, e.g., the weighted $\alpha$-fair policies of \cite{MoWa00}.  Only recently, authors have begun to consider these policies in the context of switched networks; 
to the best of our knowledge, application of $\alpha$-fairness to switched networks was first made by \cite{ShWi11} and \cite{Zh12}.  (Added in revision:  \cite{SrLi} has shown stability of a per-link scheduling policy for multihop networks.) 

\iffalse
In recent years, there has been much progress on the analysis of decentralized throughput optimal algorithms for switched networks. 
%
In particular, see \cite{JiWa10} and \cite{ShSh12}, where
%
the proofs are given for single hop networks.  When the components of networks communicate, as in the present paper, further analysis is needed. We hope that the work here will be helpful in extensions from the single hop setting to multihop networks in analyzing networks of communicating components.
\fi

\subsection{Organization} 
The remainder of the paper is organized as follows. In  Section \ref{Model}, we define a family of FIFO switched networks and the Proportional Scheduler. The
main result of the paper, Theorem \ref{mainthrm}, states that the
corresponding network is positive recurrent for all subcritical arrival rates. In Section \ref{Comparison}, we discuss the properties of the Proportional Scheduler in comparison to BackPressure; the
section is not needed to understand Theorem \ref{mainthrm}, but it is important in order to understand its consequences.  
In Section \ref{PROOF}, we prove Theorem \ref{mainthrm}.  We begin by characterizing
the fluid model that is associated with the Proportional Scheduler and then,  in (\ref{def:L.fun}-\ref{def:entropy}) and Proposition \ref{HDiff}, define a Lyapunov function for the fluid model and calculate its derivative.    This is applied, in Theorem \ref{FluidStable}, to prove fluid stability,  from which we conclude, in Proposition \ref{propFMQN}, that the corresponding
stochastic network is positive recurrent.  Proposition \ref{propFMQN} and certain other steps in the proof of Theorem \ref{mainthrm} will be proved in the appendix.

\section{FIFO Network Model, Proportional Scheduler, and Main Result}\label{Model}

\subsection{Network and Scheduling Set Notation}
\label{subsection2.1}
Let $\mJ$ be a finite set of queues, with cardinality $|\mJ|$ and  indexed by $j$.
A \emph{schedule} is a vector $\sigma=(\sigma_j : j\in\mJ) \in \bZ_+^{|\mJ|}$, where $\bZ_+^{|\mJ|}$ denotes the non-negative integers.  We will denote by
$\mS$ a finite set of schedules satisfying: (1) If $\vecsigma\in\mS$, then
$\tilde{\vecsigma}\in\mS$ for each $\tilde{\vecsigma}\leq\vecsigma$   
(with vector inequalities  $\tilde{\vecsigma}\leq\vecsigma$ being interpreted componentwise, i.e., $\tilde{\sigma}_j \leq \sigma_j$ for all $j\in\mJ$).  Note that
the vector of all zeros belongs to $\mS$.
(2) For each $j\in\mJ$, there exists some schedule $\vecsigma\in\mS$ such that $\sigma_j>0$.    For a vector $\vecQ=(Q_j: j\in\mJ) \in \bZ_+^\mJ$, we denote by $\mS_Q$ the schedules in $\mS$ with $\sigma_j \leq Q_j$ for $j\in\mJ$, and 
denote by $\sigma_{\max}:=\max\{ \sigma_j : j\in\mJ, \vecsigma\in\mS \}$
the \emph{maximum component} in the set of schedules.
We define $\coS$ to be the convex combination of points in $\mS$ and assume that $\coS$ has non-empty interior.
The \emph{subcritical region} $\mC$ of this network is the interior of $\coS$.
(In the setting of multiclass queueing networks with a fixed service rate at each queue, $\coS$ becomes a rectangle with faces parallel to the axes.)  

Each packet in the network is assumed to belong to a class at a given time $t$. The notation of a packet\rq{}s class will be used to uniquely identify the route of the packet, the queue it is at, and the stage along its route.
We denote by $\mK$ the set of classes of packets in the network.
A route through the network is a vector of classes
$r=(k^r_i : i=1,\ldots,|r|)\in\mK^{|r|}$, with size $|r|\in\bN$. Each class is assumed to occur along a unique route, with the class occurring  exactly once along its route; $\mR$ denotes the set of routes through this network.

With each class $k\in\mK$, we associate a unique route denoted $r(k)$ and a unique queue $j(k)$. For notational convenience, we add an additional ``outside\rq{}\rq{} class denoted by $\out$:  For each class $k\in\mK$, we let the function $b(k) \in\mK \cup \{\out\}$ denote the class \emph{before} class $k$ on route $r$; if $k$ is the first class on a route, we then set $b(k)=\out$. Similarly, the function $n(k) \in\mK \cup \{\out\}$ will denote the \emph{next} class on route $r$. If $k\in\mK$ is the last class on a route, then $n(k)=\out$. 
%In this way, we denote the class that a packet belongs to before and after service as a class $k$ packet at queue $j(k)$.

For a given route $r\in\mR$, we define the input class $i(r)\in\mK$ to be the first class on route $r$, i.e., $i(r)=k^r_1$, and the output class $o(r)\in\mK$ to be the last class on route $r$, i.e., $o(r)=k^r_{|r|}$. The subsets of input and output classes are denoted by $\mK^i$ and  $\mK^o$, respectively.
For notational convenience, we write $k \in j$ to indicate that class $k$ is at queue $j$, i.e., $j(k)=j$, and $k\in r$ to indicate that class $k$ occurs on route $r$, i.e., $r(k)=r$. Unless stated otherwise,  $| \cdot |$ denotes the $L^1$ norm.

\iffalse
\begin{remark}\label{ScheduleRemark}
\NW{In the model defined above we assume, for concreteness, that a schedule is deterministic. However, we note that a straight-forward extension would allow random service. Such a setting would be of interest for wireless communication networks, where packets transmitted by the sender may fail to transmit due to interference. Here each schedule would be a bounded random variable on $\bZ_+^{|\mJ|}$ with mean $\sigma$. We would then let the set $\mS$ index the set of mean values for these random variables. Our positive recurrence result, Theorem 1, holds as since the functional law of large numbers and thus our fluid model and fluid stability results which apply in this case.}
\end{remark}
\fi

\iffalse
If $x\in\mJ$, it denotes the amount of packets arrived at queue $x$,  in case $x\in\mR$, it gives  the amount of packets arrived at route $x$ and, finally, if $x\in\mK$, it is the amount of arrived packets of class $x$.
In general to easy the notation for the previous three cases we are going to use the indexes $j\in\mJ$, $r\in\mR$ and $k\in\mK$, having $A_j(t)$,  $A_r(t)$ and $A_k(t)$, respectively.
\fi

\subsection{Network Quantities and Equations}\label{FIFOeq}

Our principal objects of interest are discrete time FIFO networks and the
Proportional Scheduler, which are described here and in the next subsection.  We begin by introducing the state primitives for first-in, first-out (FIFO) networks having time index $t\in\bZ_+$.  Analogous primitives and equations will be employed in Section
\ref{PROOF} for the corresponding fluid model, where time will instead be continuous and the model will be both continuous and deterministic.  

%We now define the state primitives for a first-in, first-out (FIFO) queueing network.  Analogous primitives and equations will be employed in Section
%\ref{PROOF} for the corresponding fluid model.   In  
%(\ref{eq:Incr}-\ref{eq:ADk}) below, time is to be interpreted as discrete; in
%Section \ref{PROOF}, the time index will be continuous.

Throughout the paper, the indices $j, k,r$ will be used to refer to, respectively, queues $\mJ$, classes $\mK$ and routes $\mR$. Regardless of the index $x=j,k,r$, we will denote by $A_x(t)$ the cumulative number of arrivals by time $t$, by $D_x(t)$ the cumulative departures by time $t$, and by $Q_x(t)$ the queue size at time $t$.
For example, $A_j(t)$ is the total number of arrivals at queue $j$ by time $t$, $D_k(t)$ is the total number of packets that have departed from class $k$, and $Q_r(t)$ is the number of packets that have arrived at but not departed from route $r$.
%In addition, we must give the ordering of packets within a FIFO queue. 
%We do this by recording the order that packets of each class have been (or will be) served within each queue after $s$ units of service.
For each  $k\in\mK$ and $s\in\bN$, the function $\Gamma_k(s)$ denotes the number of packets of class $k$ that will be served after $s$ packets are served from queue $j(k)$.

The processes $A_x(t), D_x(t)$\text{ and }$\Gamma_k(t)$ are non-negative and non-decreasing, with $A_x(0) = D_x(0)=\Gamma_k(0) =0$ for $x\in \mJ\cup \mK\cup \mR$ and $k\in\mK$.  The queue size process $Q_x(t)$ is non-negative for $x\in \mJ\cup \mK\cup \mR$. Given these natural conditions, the following fundamental equations hold for a FIFO switched network.
\begin{align}
& Q_x(t) = Q_x(0)+ A_x(t) - D_x(t), \label{eq:Incr}\\
 %\label{eq:disc.Qj}
&\sum_{k \in j} \Gamma_k(t) = t, \label{eq:gammak} \\
\Big( &\frac{D_j(t) -D_j(s)}{t-s} : {j\in\mJ}\Big) \in \coS   \label{eq:Dcov} \quad \text{for } t>s, \\
&D_k(t) = \Gamma_k(D_j(t)), \label{eq:Dk}\\%\label{eq:Ak}
&A_k(t) + Q_k(0) = \Gamma_k(A_{j(k)}(t) + Q_{j(k)}(0)), \label{eq:Ak}\\ %
&A_k(t) = D_{b(k)}(t), \label{eq:ADk}
\end{align}
where $x\in \mJ\cup \mK\cup \mR$, $j\in\mJ$, $k\in\mK$ and $r\in\mR$ for (\ref{eq:Incr}-\ref{eq:Ak}), and $k\in\mK \backslash \mK^i$ in \eqref{eq:ADk}.
Finally, we define the arrivals/departures for routes by 
\begin{align}\label{eq:ADr}
A_r(t) &= A_{i(r)}(t), \qquad\qquad D_r(t) = D_{o(r)}(t).
\end{align}

The above equations (\ref{eq:Incr}-\ref{eq:ADk}) can be interpreted as follows: (\ref{eq:Incr}) is standard, \eqref{eq:gammak} states that the $t$th packet served from queue $j$ is from class $k\in j$, and \eqref{eq:Dcov} states that the departure process must be achievable within the constraints of the network scheduling. 
In \eqref{eq:Dk}, $\Gamma_k(D_j(t))$ is the number of class $k$ packets served after $D_j(t)$ units of service;
% which by definition is exactly $D_k(t)$. 
\eqref{eq:Ak} states that all of the packets that were originally at class $k$ or arrived there by time $t$ will have been served after $A_j(t) + Q_j(0)$ packets have been served at  $j$ and so, if the $A_j(t)$th packet arrival is of class $k$, then the $(A_j(t) + Q_j(0))$th packet departure is also of class $k$, i.e., the queueing discipline is FIFO; \eqref{eq:ADk} gives the routing between classes.
Further relationships can also be deduced from the above equations. For instance, \eqref{eq:gammak} and \eqref{eq:Ak} imply
%\begin{align*}
% & A_j(t) = \sum_{k \in j} A_k(t),&
%  D_j(t) &= \sum_{k \in j} D_k(t), %\label{eq:Dj}
%\end{align*}
$A_j(t) = \sum_{k \in j} A_k(t)$ and 
 $ D_j(t) = \sum_{k \in j} D_k(t)$, %\label{eq:Dj}
and \eqref{eq:Incr}, \eqref{eq:ADr} and \eqref{eq:ADk} imply
%\begin{align*}
%Q_r(t) = \sum_{k\in r} Q_k(t).
%\end{align*}
%
$Q_r(t) = \sum_{k\in r} Q_k(t)$.

We remark that, when the service rates $\sigma_j$ are constant, the FIFO property given by \eqref{eq:Ak} is equivalent to the FIFO property (2.5) in \cite{Br96a}. The FIFO condition there is given in terms of the workload at the queue; such an interpretation of workload does not immediately transfer to our setting, since service rates can vary.   Also note that (\ref{eq:Ak})
does not restrict the order in which packets that have arrived simultaneously at a queue are served; we shall allow any such order.

%This motivated the modification and extension to the definition above.

%Time here is slotted, with $t\in\bZ_+$. 
%For any discrete time process $X(t)$, we apply the notation $\Delta X(t) := X(t)-X(t-1)$. 
In addition to equations (\ref{eq:Incr}-\ref{eq:ADk}), we also assume that the number of arrivals at each route $r$ over different times is i.i.d. with mean $a_r\in (0,\infty)$, and that the arrivals at different routes occur independently.  
%
%\begin{align}\label{Aeq:Arrive}
%\{ \Delta A_r(t) : t\in\bZ_+ \} &\;\;\text{are sequences iidrvs with expectation } a_r\text{ for }r\in\mR
%\end{align}
%
For each class $k\in r$, we set $a_k = a_{r}$, and  for each queue $j\in\mJ$, we set $a_j = \sum_{k \in j} a_k$, with $\veca=(a_j : j\in\mJ)$.
We assume that the time required for the service of each packet is
deterministic and equal to $1$.  
With this in mind, we say that the mean arrival vector $\veca$ is 
\emph{subcritical} for a given network when $\veca \in \mC$;
observe that, since $\mC$ is open, this implies that $(1+\epsilon)\veca \in \mC$
for some $\epsilon > 0$. This provides a natural extension to the definition of subcriticality used for multiclass queueing networks (with a fixed allocation of service for each queue).

%
%\begin{equation*}
%a_j = \sum_{k \in j} a_k.
%\end{equation*}

%\begin{remark}
%The FIFO queueing network equations, above, do not uniquely specify a Markov process for a discrete time queueing network. To do this we must fully define the initial state of the system by specifying $\Gamma_k(s)$ for $0\leq s \leq Q_j(0)$. Further, we must define a tie-breaking rule. Here, when two or more packets from different classes arrive at a queue then we specify in the functions $\Gamma_k$  which arrive occurred first.
%\end{remark}
%
We note that for a network with the FIFO property, one needs to
specify the initial state of the system by including $\Gamma_k(s)$, for $0\leq s \leq Q_j(0)$, in order to uniquely specify the corresponding Markov process.  
As mentioned above, one also needs to define a tie-breaking rule when two or more packets from different classes arrive at a queue at the same time, 
for which we specify in the functions $\Gamma_k$  which packet arrived  ``first".

\subsection{Proportional Scheduler}
%We describe how a schedule is chosen at a given time. 
A scheduling policy is a sequence of schedules $\vecpi(t)\in\mS$ with $\vecpi(t) \leq \vecQ(t)$ component-wise
%for $j\in\mJ$ and $t\in\bZ_+$
 that determines the service of packets at each queue, i.e., $\vecD(t)-\vecD(t-1)=\vecpi(t)$. The Proportional Scheduler, which is the main focus of this paper, is defined as follows.
 %As we discussed previously, the most prevalent scheduling policies for multihop switched networks is the BackPressure policy. With some additional notation, we define this policy and then compare a number of distinctions between these two policies. The most important feature noted is that the Proportional fair policy can be implemented in FIFO networks with a per-link queueing opposed to BackPressure which must be implemented with per-route queueing.
%
%
For $\vecQ= (Q_j: j\in\mJ)\in\bZ_+^{|\mJ|}$, let $\vecsigma(\vecQ)= (\sigma_j(\vecQ) : j\in\mJ)\in \coS$ be a solution to the following optimization problem:
\begin{align}\label{PFOpt}
&\text{maximize}\qquad\sum_{j\in\mJ} Q_j \log \sigma_j\qquad
\text{over} \qquad  \vecsigma \in <\! \mS_{\vecQ} \! >;
\end{align}
when $Q_j=0$, set $\sigma_j(\vecQ)=0$, with the convention that
$0\log 0 = 0$. (Note that $\vecQ \mapsto \vecsigma(\vecQ)$ is invariant under scalar multiplication, i.e., $\vecsigma_j(c \, \vecQ) = \vecsigma_j(\vecQ)$ for $c>0$.) 
The solution to this optimization need not belong to the set of schedules $\mS$. However, since $\vecsigma \in <\! \mS_{\vecQ} \! >$, $\vecsigma$ can be expressed as a convex combination of points in $\mS_{\vecQ}$, i.e., there exists random $ \boldsymbol{\pi}(\vecQ)= ( \pi_j(\vecQ) : j\in\mJ)$, with support in $\mS$ and such that, for $j\in\mJ$,
\begin{subequations}\label{DServe}
\begin{equation}
\bE \pi_j(\vecQ) = \sigma_j(\vecQ). \label{eq:Dep}
\end{equation}
The Proportional Scheduler is defined to be any policy having the sequence of
schedules $\vecpi (t)$, where
\begin{equation}\label{Dchange}
 \vecsigma(t) = \argmax \quad \sum_{j\in\mJ} Q_j(t-1) \log \sigma_j \quad\text{over} \quad \vecsigma \in <\! \mS_{\vecQ(t-1)} \! >
\end{equation}
\end{subequations}
and $\vecpi (t)$ satisfies the analog of (\ref{eq:Dep}).
(The only sources of randomness in the paper are the choice of $\vecpi$ in (\ref{eq:Dep}) and the random exogenenous arrivals $A_r(t)$ in (\ref{eq:ADr}).)

Note that $\vecsigma (t)$ is uniquely defined because of the strict concavity of the objective function in (\ref{PFOpt}) and because $\sigma_j(t)=0$ when $Q_j(t)=0$, although $\vecpi (t)$ need not be uniquely defined (depending on $\mS$).
Also note that the Proportional Scheduler (and MaxWeight) policies reduce to the FIFO policy of a multiclass queueing network, having a fixed service rate at each queue, when $\coS$ is a rectangle with faces parallel to the axes.

We will refer to any discrete time Markov chain satisfying the FIFO switched network equations  (\ref{eq:Incr}-\ref{eq:ADk}), along with the policy in 
% \eqref{Aeq:Arrive} and 
\eqref{DServe} and the arrival and service assumptions in the next to last paragraph of  Subsection \ref{FIFOeq},
as a \emph{proportional switched network}.  The state space of the Markov chain is assumed to be that induced by the number of packets at each queue together with their respective ordering within the queue.  (With a slight abuse of terminology, we blur here the distinction between the Markov chain and the underlying switched network.)
Recall that positive recurrent means that there exists a state (here, the state with no packets) that is positive recurrent and that this state is visited with probability $1$ starting from any other state.
(This is the reduction to the countable state space setting of positive Harris recurrence from the general state space setting.)

\subsection{Main Result}
A standard fact is that, under any policy, a queueing network cannot be positive recurrent when the vector of arrival rates $\veca$ lies outside the network's subcritical region $\mC$.  
It follows that
%For instance, see proofs \cite[Theorem 3.1]{TaEp92} for switched networks, \cite[page 1060]{KeWi04} for Massouli\'e networks and  \cite[Proposition 4.11]{Br08}  for head-of-the-line queueing networks.  
%
$\mC$ is the network's greatest possible stability region.
% namely, to be positive recurrent for any arrival rate within the subcritical region.  
%(Or, equivalently, only a network with subcritical arrival rates
%can be positive recurrent.)
%
The main result in this paper is the following converse, which shows that the Proportional Scheduler achieves its greatest possible stability region.

\begin{theorem}\label{mainthrm}
Suppose that the vector of arrival rates $\veca = (a_j : j\in\mJ) \in \mC$.  Then the corresponding proportional switched network is positive recurrent.
\end{theorem}

Although this paper considers only proportional switched networks in the discrete time context, we remark that the analogous model can be defined in continuous time,
with packets instead requiring an exponential amount of service and the policy being updated immediately after each change in the state of the network.  The proof of 
Theorem \ref{mainthrm} relies on the application of fluid models to the proportional switched network and is insensitive to whether time is discrete or continuous, and so will carry over to the continuous time setting.  As in the context of FIFO multiclass queueing networks, the main requirement for proportional switched networks on their service times is that their 
means be the same for different classes at the same queue, although the means may differ among different queues.  In our setting, the assumption that service times are all of unit size at different queues is due to the setting of switched networks, rather than intrinsic mathematical requirements for Theorem \ref{mainthrm}.

As is the case for proportional switched networks, the well-known BackPressure networks are also positive recurrent for all subcritical arrival rates.
In the next section, we will compare the implementability of the two models.

\section{Comparison with the BackPressure}
 \label{Comparison}\label{BP Compare}
%Here, we compare various properties of the Proportional Scheduler with those of the
%BackPressure policy.   
The BackPressure policy is currently the canonical policy for scheduling in multi-hop switched networks.
Here, we present examples to illustrate practical differences between this policy and the Proportional Scheduler. 
This section is not required for the remainder of the paper.
%
%; the material here is included in order to motivate and emphasize properties of the Proportional Scheduler.
%
%
%This requires us to introduce some amount of new notation, not required in subsequent sections. 
%For this reason we emphasize that this section is not required to prove the main result, Theorem \ref{mainthrm}. 
%

%Although we will not give a general mathematical analysis of the BackPressure policy, in this section we emphasize the differences of this policy with respect to the proportional fair through specific examples.
%
\subsection{Definition of BackPressure}
We define and briefly describe the BackPressure policy using the notation introduced in Section \ref{Model}, with the reader being referred to \cite{TaEp92} for a more detailed description of the policy and its properties. The network is defined there in terms of a directed graph, but the reader can easily check that the definition we give below is equivalent; this alternative format is employed to avoid unnecessary complications. (Note that  queues $j\in\mJ$ are referred as links in \cite{TaEp92}, and  classes are referred to as queues there.) 
 
%An important difference between the BackPressure policy and our \emph{Proportional Scheduler} is that, for the latter, queues can have  multiple buffers,  one for each class $k\in j$,
%with scheduling being assigned according to different buffer priorities. 

For a given buffer size distribution $\vecQ=(Q_k: \, k\in\mK)$, the BackPressure policy can be defined as follows:
\begin{enumerate}
\item For each queue $j \in \mJ$, introduce the weights 
\begin{equation}\label{BPwOpt}
w_j(\vecQ)= \max_{k \in j} \left\{   Q_k - Q_{n(k)}\right\},
\end{equation}
with $Q_{n(k)}= 0$ for $n(k)=\out$.  
Denote by $k^*_j(\vecQ)$ one of the classes where the above maximization is achieved.
\item Over the set of schedules $\mS$, solve the optimization
\begin{equation}\label{BPOpt}
\sigma^*(\vecQ) \in \argmax_{\sigma\in\mS} \;\; \sum_{j \in \mJ} \sigma_j w_j(\vecQ) \ .
\end{equation}
\item When $w_j(\vecQ)>0$,  schedule $\sigma_j^*(\vecQ)$ packets from class $k^*_j(\vecQ)$ at queue $j \in \mJ$ and no packets from any of the other classes at $j$ during the next time increment; 
when $w_j(\vecQ)\le 0$, do not schedule any packets at the queue.
\end{enumerate}
 
Assuming that the number of  arrivals of packets at each route  over different times is i.i.d. and the arrivals at different routes occur independently, then    
the BackPressure policy is positive recurrent whenever 
the arrival rates are subcritical.
This is shown by employing a quadratic Lyapunov function;
the BackPressure policy, in fact, maximizes the negative drift of
the Lyapunov function. 
%belong to the capacity region $\alpha\in \mC$. 
%This holds because if one takes as a Lyapunov function the sum of the squares of the vector $Q=(Q_k: \, k\in\mK)$
%\begin{equation}
%H^*(\vecQ)=\sum_{k\in\mK} Q_k^2
%\end{equation}
%then the BackPressure policy is designed to give the most negative drift to this function.
%

%\subsection{Instability of MaxWeight}
%\subsection{Queueing Structure of BackPressure}
\subsection{Complexity of the Proportional Scheduler versus that of BackPressure}
\label{seccomplexity}

\begin{figure}[h]%
\centering
\includegraphics[width=0.8\textwidth]{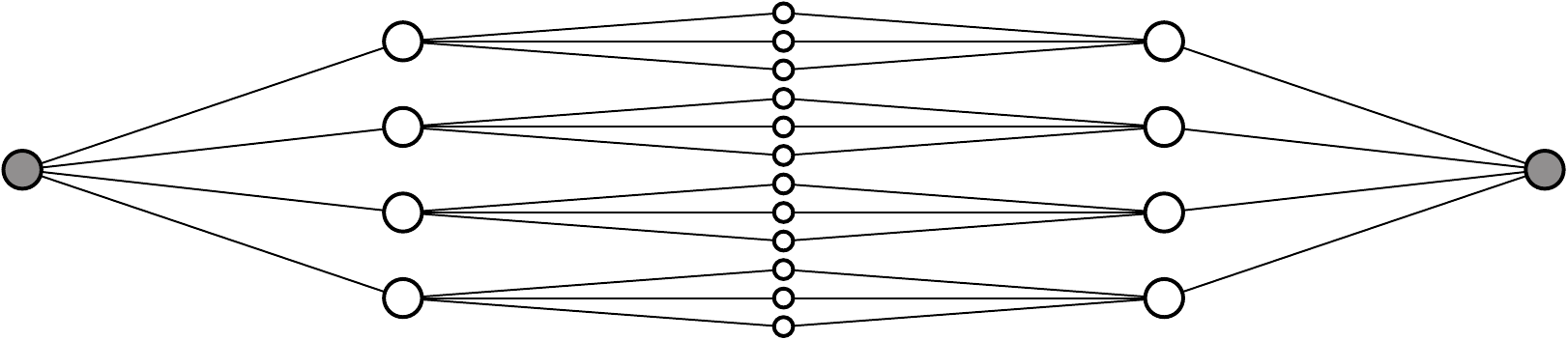}
\caption{A network of degree $d=4$ and diameter $D=4$. 
\label{Tree}}
\end{figure}%

The robust stability of the BackPressure policy is a compelling feature. However, a crucial disadvantage of the BackPressure approach is that it leads to a priority policy that requires explicit information about the classes and routes throughout the network; in many practical circumstances, the compilation of such information is not feasible.

Consider, for example, the network of diameter $D \in 2\mathbb{Z}_+$ that is constructed by joining the leaf nodes of two trees of diameter $D/2$ each, where each non-leaf node has degree $d$ -- the case of a network with degree $d=4$ and diameter $D=4$ is given in Figure \ref{Tree}. The routes in this network are the paths between the two root nodes, both from left-to-right and from right-to-left. (The two root nodes are colored grey in Figure 
\ref{Tree}.) No scheduling constraints on the nodes are assumed, and packets require unit service at each node, with each packet visiting each node along its path.  In terms of our previous
terminology, nodes here will correspond to queues, and each distinctive path emanating from a given node will correspond to a different class at that node.  (Paths that merge at a node and remain the same thereafter will correspond to the same class there.)

 For scheduling at each node,
the BackPressure policy employs computations along each of these classes. 
Focusing on one of the two root nodes for concreteness, and 
denoting by $b(d,D)$   % = |k\in j|$ 
the number of classes at the node,
  it is easy to check that
\begin{equation}\label{ExBPbound1}
b(d,D) = d(d-1)^{D/2-1} + 1 ,
\end{equation}
which grows rapidly with respect to $d$ and $D$.  If implemented in this and other settings with large networks,
the BackPressure policy will present serious memory problems. Moreover, the relevant routing information may not be available. 
%
%%For instance, queue lengths are often estimated at a road junctions for scheduling purposes; however, knowledge of car routes and thus ordering cars according to their routes is not possible.
%In a general setting, one would expect implementation problems because of the rapid growth of the number of classes needed at each site. 

In contrast to this, for the Proportional Scheduler policy, each node requires only knowledge of the number of packets destined for its $d$
adjacent nodes (and not the number of packets in each of its classes).
%
%\begin{equation*}
%b(d,D) = d \ .
%\end{equation*}
%
This amount of information is far less than that required by the BackPressure policy  in \eqref{ExBPbound1}; this quantity
does not increase as the network size increases, since it
depends only on the local structure of the network and not, for example, on the routing or
final destination of each packet. 
%
%In Figure, we plot for different values of $D$ the computational (and thus memory) requirements at the root node for BackPressure in comparison to Proportional Scheduler and MaxWeight. 

Although the network in Figure \ref{Tree} is used for reasons of exposition, 
the Proportional Scheduler is far closer in nature to the next-hop routing used in modern IP routers on the Internet, or that might be used in a wireless ad-hoc network. We note, for instance, that modern Internet routers maintain tens of FIFO queues that aggregate tens of thousands of flows (route classes) (\cite{mckeown1999islip,appenzeller2004sizing}). Furthermore, when these components are connected together to form a network, queue state information is not explicitly changed.  (We discuss this in  Subsection \ref{secdecomposition}.)

\subsection{Simulations for a Modification of the Network of Figure \ref{Tree}}
\begin{figure}[h]%
\centering
\includegraphics[width=0.70\textwidth]{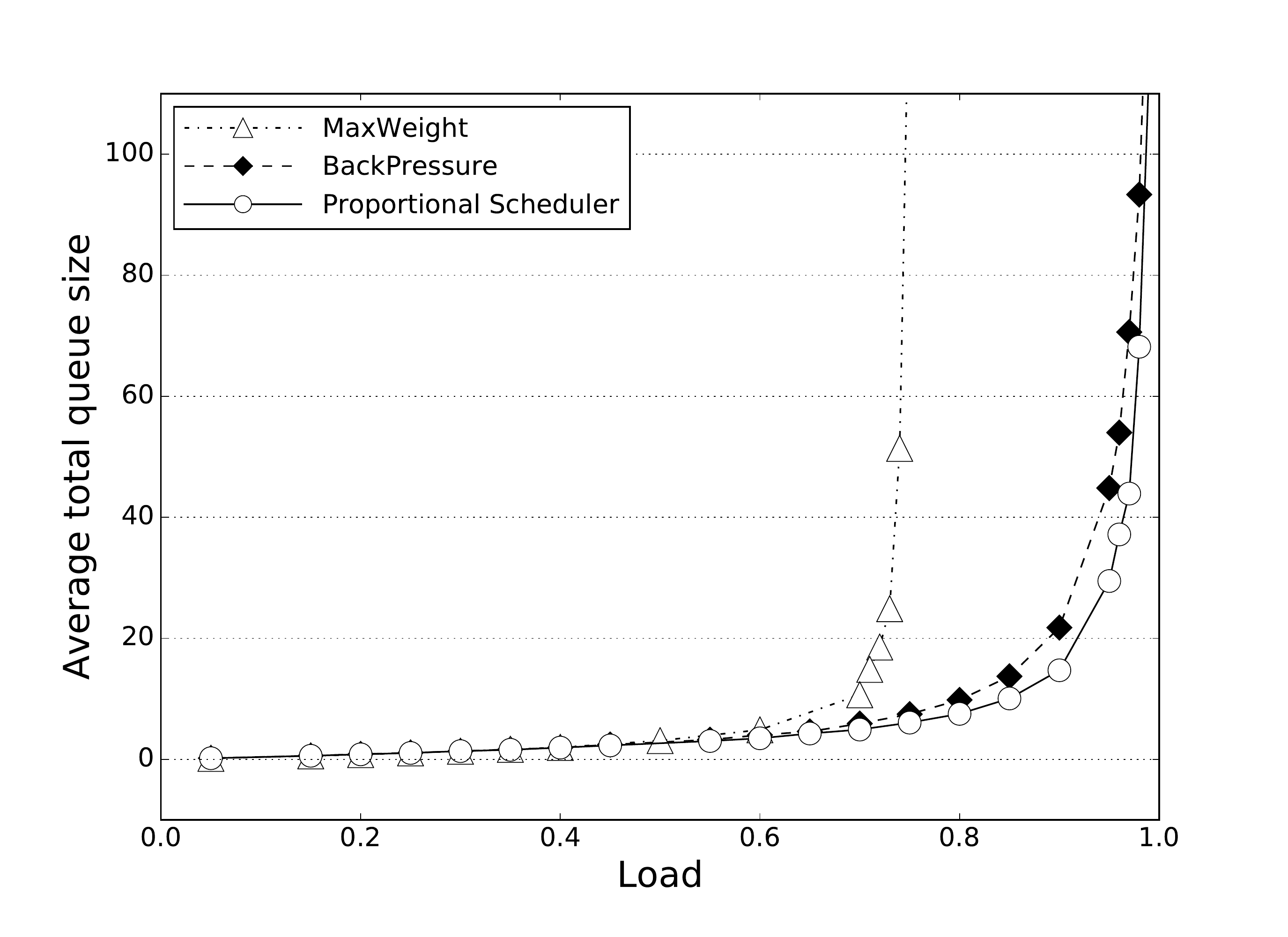}
  \caption{Total queue size of MaxWeight, BackPressure, Proportional Scheduler under a range of loads.%, Figure \ref{Tree}. 
 \label{MWSim2}}
\end{figure}

In Figure \ref{MWSim2}, we provide simulations for the average total queue size for a modification of the network in Figure \ref{Tree}.  Specifically, rather than immediately exiting from the root node at the end of its route, a packet returns to this node a geometric number of times, each time with probability $1-\delta$, $\delta \in (0,1]$, before finally exiting from the network.   (We could instead stipulate that the packet returns to the node a fixed number of times before exiting.)

For this modified network, the load at a root node is $a (1+ 1/\delta)$ if $a$ is the external arrival rate at each root node, and so $a(1+1/\delta) < 1$ is an obvious necessary condition for 
%
%\begin{equation}\label{StableCond4}
%%
%%a\Big(1-\frac{1}{\delta}\Big) < 1
%%
%
%\end{equation}
%
positive recurrence of the network, whereas, for both BackPressure and the Proportional Scheduler, the condition is sufficient by \cite{TaEp92} and Theorem \ref{mainthrm}. For parameters $(D,d,\delta)=(4,4,0.4142)$, we simulated 
in Figure \ref{MWSim2}  three policies, BackPressure, the Proportional Scheduler, and 
MaxWeight, over a range of loads. Our simulations indicate that, as expected, BackPressure and the Proportional Scheduler are maximally stable, whereas, consistent with the examples of input-queued switches in \cite{andrews2003achieving},  MaxWeight is not.

We note that the MaxWeight has a stability region that is approximately 75\% of the total capacity. ($\delta = .4142$, which is about $\sqrt{2} -1$, was chosen so that value would not be close to $1$.)  In fact, MaxWeight will lose up to 50\% of the capacity for networks of this type; as mentioned in the introduction, more detail on MaxWeight will be provided in a future work. We also note that, under all loads measured in Figure \ref{MWSim2}, the average total queue size for the Proportional Scheduler was found to be somewhat smaller than that for BackPressure. We will discuss this relationship for another family of networks in Section \ref{delay}.

\subsection{Delay on Long Routes}\label{delay}
As illustrated in Subsection \ref{seccomplexity},
the complexity of the BackPressure policy is an issue when there are many routes.  Other complications might also arise, even when all packets have the same route.
Consider, for example, the network in Figure \ref{Line}, which
%
%\footnote{Switching in this linear network is not to be confused with the simultaneous resource possession in linear networks considered in the context of Bandwidth Sharing networks by \cite{BoPr03}.} 
%
is a special case  of the linear networks considered by \cite{BSS11} (in Theorem 2) and \cite{St11}, for the BackPressure policy. 
%
%(Switching in this linear network should not be confused, however, with the simultaneous resource possession in linear networks considered in the context of Bandwidth Sharing networks by \cite{BoPr03}.) 
%
\iffalse
\begin{figure}[h!]
  \centering
\includegraphics[width=0.95\textwidth]{Spin_Line.pdf}
\caption{A linear network with $J$ links. Each packet must pass through each link.
\label{Line} }
\end{figure}
\fi

\iffalse
\begin{figure}[h]%
\centering
\begin{minipage}[b]{0.45\textwidth}
\includegraphics[trim=0 -25 0 0, clip, width=1.1\textwidth]{Spin_Line.pdf}
\vspace*{1,5cm}
\caption{A linear network with $J$ links. Each packet must pass through each link.
\label{Line} }
 \raisebox{1.5ex}{\kern-4.45cm \tiny \bf (on the right)}%
\label{BP2PS-Line}
%\vspace*{0.25cm}
\end{minipage}%
\qquad
\begin{minipage}[b]{0.5\textwidth}
 \vspace*{-1cm}
  \includegraphics[width=1.1\textwidth]{BP2PS-tandem-net-Bernoulli+Poisson.pdf}
 \caption{Queue lengths for the linear network with $J=20$ links.}%
 \vspace*{0.42cm}
\end{minipage}%
\end{figure}%
\fi

\begin{figure}[h]%
\centering
\begin{minipage}[b]{0.45\textwidth}
\includegraphics[trim=0 -25 0 0, clip, width=1.1\textwidth]{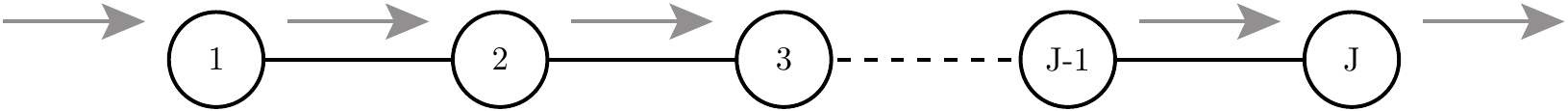}
\caption{(Above) A linear network with $J$ links. Each packet must pass through each link.
\label{Line} }
\vspace*{1,00 cm}
 \caption{(Right) Queue lengths for the linear network with $J=20$ links.}%
 \raisebox{1.5ex}{\kern-4.45cm \tiny \bf (on the right)}%
\label{BP2PS-Line}
\vspace*{0.25cm}
\end{minipage}%
\qquad
\begin{minipage}[b]{0.5\textwidth}
  \includegraphics[width=1.1\textwidth]{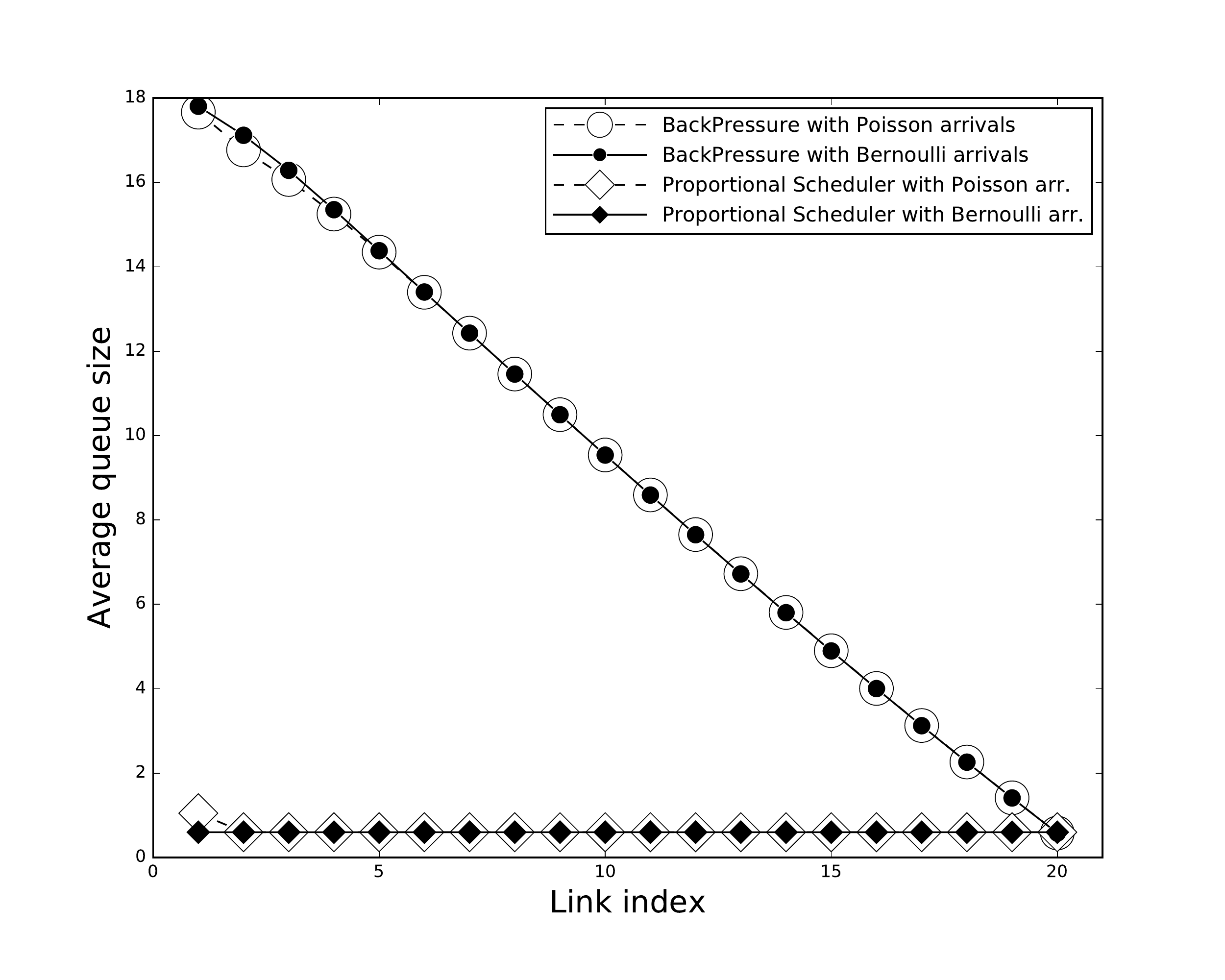}
 \vspace*{-1cm}
\end{minipage}%
\end{figure}%

In this example, there are $J\ge 2$ links in series, with packets entering at link $1$ and being sequentially processed through links $1,2,\ldots,J$;
there are no constraints preventing simultaneous service at different links and,
once served, a packet moves to the next link along its route.
\iffalse
We assume there is a maximum number of allowed failures at a given link, say $\bar q$.%
\footnote{We remark that, in this setting, we can define BackPressure policy with the probability $1-q_j$ in expression \eqref{BPOpt}. 
Further, with the assumed technical condition, the probabilistic routing in this example is a special case of the routing structure considered in this paper for the proportionally scheduler.}
\fi
%A practical interpretation of such a network could be as a series of wireless routers with a given  probability of failing to sending packets. 
Arrivals at queue $1$ occur according to a Bernoulli process, with mean parameter $a>0$, and one packet is served per unit time at each nonempty queue. %(Recall that a Bernoulli process is a counting process whose increments are i.i.d. Bernoulli random variables.)
%We assume that each link is subcritically loaded, i.e., that $a  <1$ in this case. 
\cite{BSS11} and  \cite{St11} showed that,
 if the network is subcritically loaded with $1/2 <a<1$,  then, under the BackPressure policy,
 the sum of the expected queue sizes in equilibrium grows quadratically in $J$, that is,
\begin{equation}
\label{eqofprop1}
\sum_{j=1}^J \bE \big[ Q_{j} \big] \geq  c \, J^2 
\end{equation}
for some constant $c>0$ not depending on $J$.
\iffalse
\begin{proposition}\label{pr:linear.network.BP}
In a linear network with $J\ge 2$ links, if the network is subcritically loaded and $\alpha_j > 1/2$ for all $j$, then, under the BackPressure policy,
 the sum of the expected queue sizes in equilibrium grows quadratically in $J$; in particular,
\begin{equation}
\label{eqofprop1}
\sum_{j=1}^J \bE \big[ Q_{j} \big] \geq  c \, J^2 
\end{equation}
for some constant $c>0$ not depending on $J$.
\end{proposition}
\fi
Since one unit of time is required to serve each packet at a queue, it is immediate from (\ref{eqofprop1}) that the flow delay for packets satisfies the same lower bound.

%%The basic idea behind the demonstration of (\ref{eqofprop1}) is that  BackPressure only serves a packet on the $j$th link when $Q_j  - Q_{j+1} > 0$, which will occur in equilibrium with probability $a$ (since the average arrival and departure rates are equal), whereas, in equilibrium, $Q_j  - Q_{j+1} \ge -1$ must always hold (since service from queue $j$ occurs only when $Q_j - Q_{j+1} > 0$).  Since $a > 1/2$, there will be linear growth in the expected size of successive queues from the last queue $Q_J$ to the first queue $Q_1$, which implies the quadratic growth in $J$ of the sum of the expected queue sizes.
%%
%%and, in order to be able to
%%process all packets on the link,  $Q_{j}-Q_{j+1} > 0 $ needs to occur with positive probability. 
%%%Therefore, we expect each queue to be a constant bigger than the previous queue.
%%We therefore expect queue sizes along the route to increase linearly from the last queue $Q_J$ to the first queue $Q_1$, which implies the 
%%quadratic growth in terms of $J$ of the total queue size and delay.
%%%
%\begin{BD}Generally\end{BD}, when routing is more involved,
%one should expect the sum of the queue sizes in equilibrium to grow quadratically with route length for routes with close to critical arrival rates since, as before,  
%BackPressure will not serve links with a negative queue size differential.
%%In informal terms, the BackPressure policy chooses not to serve these links as it assumes that providing service to them inhibits the progress of packets at different parts of the network.
%(For details, see \cite{BSS11} or \cite{St11}.)

%

The basic idea behind the demonstration of (\ref{eqofprop1}) is that 
BackPressure only serves a packet on the $j$th link when $Q_j  - Q_{j+1} > 0$, which will occur in equilibrium with probability $a$ (since the average arrival and departure rates are equal), whereas, in equilibrium, $Q_j  - Q_{j+1} \ge -1$ must always hold (since service from queue $j$ occurs only when $Q_j - Q_{j+1} > 0$).  Since $a > 1/2$, there will be linear growth in the expected size of successive queues from the last queue $Q_J$ to the first queue $Q_1$, which
implies the quadratic growth in $J$ of the sum of the expected queue sizes.
%
%and, in order to be able to
%process all packets on the link,  $Q_{j}-Q_{j+1} > 0 $ needs to occur with positive probability. 
%%Therefore, we expect each queue to be a constant bigger than the previous queue.
%We therefore expect queue sizes along the route to increase linearly from the last queue $Q_J$ to the first queue $Q_1$, which implies the 
%quadratic growth in terms of $J$ of the total queue size and delay.
%%
In general, when routing is more involved,
one should expect the sum of the queue sizes in equilibrium to grow quadratically with route length for routes with close to critical arrival rates since, as before,  
BackPressure will not serve links with a negative queue size differential.
%In informal terms, the BackPressure policy chooses not to serve these links as it assumes that providing service to them inhibits the progress of packets at different parts of the network.
(For details, see \cite{BSS11} or \cite{St11}.)

The Proportional Scheduler  (as well as any other work conserving scheduler)  will exhibit completely different behavior for this linear network since there are no constraints preventing simultaneous service at different links. In equilibrium, queue $1$ will have either $1$ or $0$ packets corresponding to whether or not an arrival occurred in the last time slot. Therefore, the output of a queue is a Bernoulli process that is independent of the current state of the queue. Arguing inductively, this implies that the queue sizes at different queues are independent, and hence that the sum of the queue sizes for the network is binomially distributed with parameters $J$ and $a$. Consequently, under the Proportional Scheduler policy, the sum of the expected value of the queue sizes in equilibrium grows only linearly in $J$, with
\begin{equation}\label{eqofprop2}
\sum_{j=1}^J \bE \big[Q_j \big]   = a J .
\end{equation}
%The bound in \eqref{eqofprop1} is significantly worse than that in (\ref{eqofprop2}),
%which illustrates the inefficiency of the BackPressure policy for even this elementary
%network.
%as well as any other work-conserving policy,
\iffalse
\begin{proposition}\label{pr:linear.network.PF}
In a linear network with $J$ links, if each link is subcritically loaded, then, under the Proportional Scheduler policy, the total queue size of the network in equilibrium grows only linearly with $J$, i.e.,
\begin{equation*}
\sum_{j=1}^J \bE \big[Q_j \big]   \leq c \, J 
\end{equation*}
for some constant $c>0$ not depending on $J$.
\end{proposition}
\fi
%The delay properties are trivialized by independence between queues achieved under the Proportional Scheduler. 
%
%See \cite{ShWaZh12}, for a non-trival instance which exploits similar structure to obtain improved delay bounds.
%
The above reasoning that was employed for \eqref{eqofprop2} is an elementary variant of Burke's Output Theorem (as in, e.g., \cite{Bu56} or \cite{HsBu76}). 

In Figure \ref{BP2PS-Line}, we have simulated the linear network in Figure \ref{Line}, with $J=20$ queues, and assuming an arrival rate of 0.6 of packets into the network and a maximal service rate of 1 at each node.  In addition to the case of Bernoulli input discussed above, Poisson input is also simulated for both BackPressure and the Proportional Scheduler. 
% (100 trials, each of time length $10,000$ and starting from the empty state, were employed; only a much shorter time scale is needed for the network to approach equilibrium.) 
As expected from the above discussion, for BackPressure, the average queue size is proportional to the number of hops from the terminal destination while, for the Proportional Scheduler, the average queue size is constant; in neither of these cases, is the average queue size sensitive to the distribution of the input except for the first few links, as illustrated by the almost perfect overlay of the pairs of graphs.
It follows that, as the length of the linear network grows,  queueing delay grows quadratically for BackPressure and linearly for the Proportional Scheduler.  (Similar observations were made in the simulations conducted by \cite{BSS11}.)

\subsection{Decomposition}
\label{secdecomposition}
Decomposition is essential to the decentralized implementation of a policy. 
A property of BackPressure optimization is that, once queue length comparisons have been made between links (see \eqref{BPwOpt}), the optimization can be decomposed. In particular, if the scheduling of one subset of queues $\mJ_1$ does not effect the scheduling choice of a complementary subset $\mJ_2$ (i.e., $\mS=\mS_1\times\mS_2$), then the BackPressure optimization can be decomposed as 
\begin{equation*}
\max_{\sigma\in\mS} \Big\{ \sum_{j \in \mJ} \sigma_j w_j(\vecQ) \Big\}= \max_{\sigma\in\mS_1} \Big\{  \sum_{j \in \mJ_1} \sigma_j w_j(\vecQ) \Big\} + \max_{\sigma\in\mS_2} \Big\{  \sum_{j \in \mJ_2} \sigma_j w_j(\vecQ) \Big\}.
\end{equation*}
The sub-problems involving the two optimizations on the right can then be solved independently, leading to a decomposed implementation of the policy. We remark, however, that the above optimization may not completely decompose since the weight calculations $w_j(Q)$ typically require comparisons with the sizes of upstream queues, and hence some information exchange will occur between network components. 

No such queue size comparisons are required for the proportional fair optimization that is employed to define the Proportional Scheduler, and so the optimization can be completely decomposed as
\begin{equation*}
\max_{\sigma\in\coS} \Big\{  \sum_{j \in \mJ} Q_j \log \sigma_j \Big\}=  \max_{\sigma\in <\! \mS_1\!>}\Big\{ \sum_{j \in \mJ_1} Q_j \log \sigma_j \Big\} + \max_{\sigma\in <\! \mS_2\!>} \Big\{  \sum_{j \in \mJ_2} Q_j \log \sigma_j \Big\}
\end{equation*}
for complementary components having 
independent scheduling.
This leads to more potential applications when the network decomposes, in comparison with the BackPressure policy.

\section{Maximum Stability Proof} \label{PROOF}
In this section, we prove Theorem \ref{mainthrm}, for which we employ fluid model
techniques introduced in  \cite{Da95}.
%and then extended in \cite{Br96a}, see also \cite{Br08}. 
The basic idea is to show that the Markov chain converges under appropriate scaling to a deterministic fluid model, 
and then to show that the queue lengths of all normalized fluid model solutions converge to $0$ within
a fixed time.
This latter step requires most of the work in the proof 
and is shown by defining an appropriate Lyapunov function that
decreases to $0$ by this time.  

From a mathematical point of view, the main contribution of this argument is the extension of fluid model techniques to a case with varying service speeds.
Here, this requires identifying the Lyapunov function in \eqref{def:entropy}, showing it to be non-increasing in Corollary \ref{cor:H.non-increasing},
and then showing that this Lyapunov function in fact decreases at a fixed rate by applying Lemmas \ref{bound.change.Q} and \ref{lm:lowerbound.q.sigma}. 
The proofs of certain technical lemmas and propositions are relegated to the appendix.

\subsection{Fluid Model}\label{fluidmodel}

We will employ, in Section \ref{PROOF}, a fluid model, with the functions  
$A_x(t)$, $D_x(t)$, $Q_x(t)$, and $\Gamma_k(t)$, 
that satisfies the same conditions (\ref{eq:Incr}-\ref{eq:ADk}) as the discrete time switched network of Section \ref{Model}, and such that
$A_x(t), D_x(t)$\text{ and }$\Gamma_k(t)$ are non-negative and non-decreasing,
and $Q_x(t)$ is non-negative, with $A_x(0) = D_x(0)=\Gamma_k(0) =0$ for $x\in \mJ\cup \mK\cup \mR$ and $k\in\mK$.
Note that each of the conditions (\ref{eq:Incr}-\ref{eq:ADk}) is invariant under the ``law of large numbers" scaling, where both time and space are scaled equally. 

In addition to these conditions, the route arrival and queue departure processes of the fluid model will be assumed to satisfy
\begin{align}
A_r(t) =& a_r t \label{Fluid:1},\\
Q_j(t) >0 \quad \text{implies} \quad   D_j\rq{}&(t) = \sigma_j ( \vecQ(t))\,, \label{Fluid:3}
\end{align}
for given $a_r > 0$ and for $t>s>0$, $j\in\mJ$ and $r\in\mR$, with $\vecsigma(\vecQ)$ denoting a solution to the proportional fair optimization
\begin{equation}\label{PF}
\vecsigma(\vecQ) \in \argmax_{\vecsigma\in \coS} \,\,\sum_{j\in\mJ} Q_j \log \sigma_j.
\end{equation}
As in Section \ref{Model}, we set $\sigma_j(\vecQ)=0$ whenever $Q_j=0$ and, for each class $k\in r$ and queue $j\in\mJ$, set $a_k = a_{r}$, $a_j = \sum_{k \in j} a_k$, and $\veca=(a_j : j\in\mJ)$.
We will refer to the conditions (\ref{eq:Incr}-\ref{eq:ADk}) and (\ref{Fluid:1}-\ref{PF}) on $\vecA(t)$, $\vecD(t)$, $\vecQ(t)$, and $\vecGamma (t)$ as \emph{proportional switch fluid model equations}, or collectively, as the \emph{proportional switch fluid model} (or \emph{fluid model}, for short).

Together, the above conditions imply with a little work that, for $x\in \mJ\cup \mK\cup \mR $, the components $A_x (t)$, $D_x (t)$, $Q_x (t)$, and $\Gamma_x (t)$ are all Lipschitz continuous and thus almost everywhere differentiable.
Note that $D'_j(t)=0$ need not hold when $Q_j(t)=0$, since work may
enter and leave an empty queue of a fluid model.  Also note that a solution of the fluid model equations corresponding
to a given initial condition need not be the unique such solution.  

%The formal justification of the fluid model just described is given in Section \ref{appendFMFL} of the Appendix.

\iffalse
\begin{remark}
As in Section \ref{Model}, we employ here the convention that $\sigma_j(\vecQ)=0$ whenever $Q_j=0$.  However, there may be other optimal solutions to this optimization. In particular, it is important to note that, $D'_j(t)$ solves the proportional fair optimization \eqref{PFOpt} but, in general, $D'_j(t)=0$ need not hold when $Q_j(t)=0$. In other words, an empty queue may still be processing work in our fluid model.
\end{remark}
\fi
\subsection{Lyapunov Function and its Derivative}\label{LyaDeriv}

In this subsection, we define the entropy function $H(t)$ and prove that it has negative derivative for non-empty subcritical systems.
We first introduce the functions
\begin{align}
L(t) = 
& \sum_{\jInJ} Q_j(t) \, \log D'_j(t) \label{def:L.fun},  \\
M(t) = 
& \sum_\jInJ \sum_\kInj 
  \int_{D_j(t)}^{Q_j(0)+A_j(t)}
    \Gamma'_k(s) \log \frac{\Gamma'_k(s)}{a_{k}} \, ds \label{def:M.fun};
\end{align}
%and the \emph{entropy function} as 
%
$H(t)$ is then defined as
\begin{equation}\label{def:entropy}
H(t) 
 = L(t) + M(t) \ .
\end{equation}
The entropy function $H(t)$ builds on entropy functions in \cite{Ma07} and \cite{Br96a}.
The term L(t) is derived in \cite{Ma07} as the large deviations rate function of a reversible network that approximates proportional fairness; the term 
$M(t)$ is similar to the entropy function employed in \cite{Br96a} for FIFO 
multiclass queueing networks.

The function  $L(t)$ is well defined everywhere (with, as earlier, the convention
 that $0 \log 0 = 0$).  Note that  
 %by \eqref{Gamma.conservation},
\begin{equation}
\label{eq22'}
\begin{split}
Q_j(t)  =  Q_j(0) + A_j(t) - D_j(t) 
=  &\sum_\kInj (\Gamma_k(Q_j(0) +A_j(t)) - \Gamma_k(D_j(t))) \\ 
=  &\sum_\kInj \int_{D_j(t)}^{Q_j(0)+A_j(t)} \Gamma'_k(s) \, ds
\ ;
\end{split}
\end{equation}
the lower and upper bounds $D_j(t)$ and $Q_j(0)+A_j(t)$ for the integral in (\ref{eq22'}) can be thought of as the amounts of packet mass having departed from the queue $j$ by time $t$ and having departed by the later time when the last of the mass already at $j$ at time $t$ departs from the queue.  With (\ref{eq22'}) in mind, note that $M(t)$ is the sum over the different queues of a weighted version of the packet mass at each queue $j$ at time $t$, based on the departure rates of the mass at the queues.

The following lemma shows that the function $H(t)$ is always non-negative for any solution of the fluid model.

\begin{lemma}\label{prop:H.not.zero}
For any $t\geq0$,
$H(t)\ge 0.$
Moreover, whenever the arrival rate vector $\veca$ belongs to the interior of the scheduling set $\mC$, then $H(t)=0$ iff $\vecQ(t)=0$.
\end{lemma}
\proof{Proof}
By (\ref{eq22'}),
$$
Q_j(t) \, \log a_j =
  \sum_\kInj \int_{D_j(t)}^{Q_j(0)+A_j(t)}
    \Gamma'_k(s) \log a_j \, ds 
$$
for each $\jInJ$; adding and subtracting such terms
allows one to rewrite $H(t)$ as
\begin{equation}\label{HforPositivity}
H(t) 
 = \sum_\jInJ 
  \int_{D_j(t)}^{Q_j(0)+A_j(t)} \sum_\kInj 
    \Gamma'_k(s) \log \left(\Gamma'_k(s) \frac{a_j}{a_k}\right) \, ds
  + \sum_{\jInJ} Q_j(t) \, \log \frac{D'_j(t)}{a_j} \ .
\end{equation}
\iffalse
If we define the discrete probability distributions $\vecp$ and  $\vecq$, over the states $\kInj$, by $\vecp=(\Gamma'_k(s) : \kInj )$ and $\vecq= ( a_k/a_j : k \in j)$, we have that
$$
\sum_\kInj 
    \Gamma'_k(s) \log \left(\Gamma'_k(s) \frac{a_j}{a_k}\right)
= D(\vecp||\vecq).
$$
where $D(\vecp||\vecq)$ denotes the \emph{relative entropy} between the distributions $\vecp$ and $\vecq$.
\fi 
By  Jensen's inequality,  
\begin{equation}
\label{unlabeledeq}
\sum_\kInj 
    \Gamma'_k(s) \log \left(\Gamma'_k(s) \frac{a_j}{a_k}\right) = -\sum_\kInj 
    \Gamma'_k(s) \log \left( \frac{1}{ \Gamma'_k(s)} \frac{a_k}{a_j}\right)  \geq -\log \left(  \sum_\kInj 
    \Gamma'_k(s)    \frac{1}{ \Gamma'_k(s)} \frac{a_k}{a_j}\right) =0.
\end{equation}
From this, it follows that the first summation in \eqref{HforPositivity} is non-negative.
The non-negativity of the second term in (\ref{HforPositivity}) follows since  the service rates $D'_j(t)$ solve the optimization problem \eqref{PF} and  $\veca\in \mC$, and so, in particular,
$$
\sum_{\jInJ} Q_j(t) \, \log D'_j(t)
\geq 
\sum_{\jInJ} Q_j(t) \, \log a_j \ .
$$
Together with (\ref{unlabeledeq}), this implies $H(t)\geq 0$.

When the fluid model is subcritical, one can improve on the inequality in the last display to show that $H(t)$ is positive when $\vecQ(t) \neq 0$.  For this, note that
%
%Indeed, under the assumption of subcriticality, the arrival-rates vector belongs to the interior of the subcritical region, i.e. 
$\veca \in \mC$ implies the existence of an  $\epsilon > 0$ such that
$$ (1+\epsilon) \, \veca \in \mC \ .$$
 Since $\vecD'(t)$ is an optimal scheduling, 
it follows that, for $\vecQ(t) \not = 0$,
\begin{equation}\label{HQ}
\sum_{\jInJ} Q_j(t) \, \log \frac{D'_j(t)}{a_j}
\geq 
\sum_{\jInJ} Q_j(t) \, \log (1+\epsilon) > 0 \ .
\end{equation}
Consequently, $H(t)=0$ implies $Q(t)=0$.  Trivially, $\vecQ(t)=0$ implies $H(t)=0$.\Halmos\endproof
 
We next show that  $H(t)$ is decreasing when $\vecQ(t)\not=0$, by proving $H'(t)$ is negative.
Since we do not know in advance that $L(t)$ is sufficiently regular to
apply the Fundamental Theorem of Calculus to
$\int_s^t L'(u)\,du$, we will employ the following
%
%It is relatively straightforward to differentiate $H(t)$ and then find its intergral to be negative. However, its derivative does not exist everywhere and so we cannot \emph{a priori} assume that $H(t)$ is the integral of its derivavtive. 
%
technical lemma, which is based on Lemma 5 of \cite{Ma07}.

%[DO WE NEED T?]

\begin{lemma}\label{L.deriv.a.e.} 
There exists a time $T>0$ not depending on the initial queue state such that:

i) For $t \geq T|\vecQ(0)|$, 
%
%the function $L(t)$ admits the following derivative almost everywhere
%
\begin{equation}\label{eq:EnvelopeTheorem}
L'(t) = 
\sum_\jInJ Q'_j(t) \log D'_j(t) \quad \text{a.e.}
\end{equation}

ii) There exists a constant $\kappa_L>0$ such that, for all $t\geq s \geq 0$,  
\begin{equation}\label{eq:L.Lipschitz.up}
L(t)-L(s) \leq \kappa_L (t-s).
\end{equation}

iii) For all $t\geq s\geq T|\vecQ(0)|$,
\begin{equation*}
L(t) - L(s) \leq \int_s^t L'(u) \, du.
\end{equation*}
\end{lemma}
A proof of this lemma is given in the Appendix. 

Since the processes $\vecA(t)$ and $\vecD(t)$ are Lipschitz, it is easy to see that the term $M(t)$, which integrates an a.e. bounded function, is Lipschitz and hence a.e. differentiable.  
Thus, the following corollary is immediate from Lemma \ref{L.deriv.a.e.}.
%see that the properties of the function $L(t)$ given in Lemma \ref{L.deriv.a.e.} apply to the function $H(t)$. In particular the increments of the entropy function can be controlled above by the integral of its derivative as stated in the following corollary.

\begin{corollary}\label{H.deriv.a.e.} 
i) There exists a constant $\kappa_H>0$ such that, for $t\geq s \geq 0$,  
\begin{equation}\label{eq:H.Lipschitz.up}
H(t)-H(s) \leq \kappa_H (t-s).
\end{equation}
ii) There exists a time $T>0$ not depending on the initial state such that $H(t)$ is a.e. differentiable on $t \in [T|\vecQ(0)|, \infty)$ and, for given $s,t$ with $t\ge s\ge T|\vecQ(0)|$,
\begin{equation}\label{eq:upp.Lipschitz}
H(t) - H(s) \leq \int_s^t H'(u) \, du.
\end{equation}
\end{corollary}

 We now focus on computing $H'(t)$. The following proposition shows that  $H'(t)$ is the sum of two terms: the (unnormalized) relative entropy between route departure and arrival rates, and the (unnormalized) relative entropy between queue arrival and departure rates. This decomposition will play an important role in the proof of Theorem \ref{FluidStable}.  
%(As earlier, we set $0\log 0 = 0$; we also set $0\log (0/0) = 0$ here and note
%that the set where $\{A_j'(t) > 0, D_j'(t) = 0\}$ has measure $0$ for each $j$.
%For brevity, when no confusion is likely, we will often omit the quantifier ``a.e." in our computations.)

\begin{proposition}\label{HDiff}
For $T$ as in Corollary \ref{H.deriv.a.e.}, and $t\ge T|\vecQ(0)|$, 
\begin{equation}\label{H(t).deriv}
H'(t) = 
- \sum_\rInR  D'_r(t) \log \frac{D'_r(t)}{A_r'(t)} 
- \sum_\jInJ A'_j(t) \, \log \frac{A'_j(t)}{D'_j(t)} \quad a.e. 
\end{equation}
\end{proposition}
As earlier, we set $0\log 0 = 0$; we also set $0\log (0/0) = 0$ here and note
that the set where $\{A_j'(t) > 0, D_j'(t) = 0\}$ has measure $0$ for each $j$.
For brevity, when no confusion is likely, we will often omit the quantifier ``a.e." in our computations.
\proof{Proof of Proposition \ref{HDiff}}
Differentiating the expression for $H(t)$ in \eqref{def:entropy}
and applying \eqref{eq:EnvelopeTheorem} gives
\begin{align*}
H'(t) = 
& \sum_{j\in\mJ}\sum_{k\in j} A'_j(t) \, \Gamma'_k(A_j(t)) \log \frac{\Gamma'_k(A_j(t))}{a_k} \nonumber \\
& -\sum_{j\in\mJ}\sum_{k\in j} D'_j(t) \, \Gamma'_k(D_j(t)) \log \frac{\Gamma'_k(D_j(t))}{a_k}
+ \sum_\jInJ Q'_j(t) \, \log D'_j(t) .
\end{align*}
%
%where for the derivative of the $L(t)$ we used 
%\eqref{eq:EnvelopeTheorem} from Lemma \ref{L.deriv.a.e.} .
On the other hand, differentiation of the expressions \eqref{eq:Dk} and \eqref{eq:Ak} implies that
\begin{equation*}
D'_k(t) = D'_j(t)\Gamma'_k(D_j(t)), \qquad A'_k(t) = A'_j(t)\Gamma'_k(A_j(t)),
\end{equation*}
and substitution of these terms into $H'(t)$ gives
%that give the arrival and the departure rates of class $k$ packets to and from queue $j$, we get
\begin{align*}
H'(t) = 
& \sum_{j\in\mJ}\sum_{k\in j} A'_k(t) \log \frac{A'_k(t)}{A'_j(t) \, a_k}
-\sum_{j\in\mJ}\sum_{k\in j} D'_k(t) \log \frac{D'_k(t)}{D'_j(t) \, a_k} \nonumber \\
& + \sum_\jInJ (A'_j(t)- D'_j(t)) \, \log D'_j(t) \\
= & \sum_{j\in\mJ}\sum_{k\in j}  \left(A'_k(t) \log \frac{A'_k(t)}{a_k} 
- D'_k(t) \log \frac{D'_k(t)}{a_k} \right)
+ \sum_\jInJ A'_j(t) \, \log \frac{D'_j(t)}{A'_j(t)} \nonumber .
\end{align*}
Application of \eqref{eq:ADk} shows that the double sum above telescopes over the classes in each route;
since  $A'_r(t) \equiv a_k$ for each $k\in r$, only the terms due to 
 external departures $D'_r(t)$ do not cancel out,  which implies
 \eqref{H(t).deriv}.
\Halmos\endproof

It is now easy to see that $H'(t)\le 0$ for suitably large $t$.
\begin{corollary}\label{cor:H.non-increasing}
For $T$ as in Corollary \ref{H.deriv.a.e.} and $t>T|\vecQ(0)|$, $H'(t) \le 0$ a.e.
\end{corollary}
\proof{Proof}
For each choice of $a,\! d\geq0$, one has  $a \log (a/d) \geq a-d$ (as before,
 $0\log 0 = 0$,  $0\log (0/0)= 0$). Applying this to \eqref{H(t).deriv} implies that
\begin{equation}\label{eq:cons.flow}
H'(t) \leq  \sum_{r\in\mR} \big( A'_r(t) - D'_r(t) \big) -  \sum_{j\in\mJ} \big( A'_j(t) - D'_j(t) \big)
=\sum_\rInR Q'_r(t)
- \sum_\jInJ Q'_j(t) = 0 \ ,
\end{equation}
with the final equality following since 
%
%follows by the conservation flows, that is, 
%
the rates of change in total queue size summed over all queues and summed over all routes are equal. \Halmos\endproof

\subsection{Fluid Stability and Positive Recurrence}\label{PosRec}

We wish to show that, in an appropriate averaged sense, $H(t)$ is always decreasing at a rate bounded away from zero when $\vecQ(t)\not=0$.  It will follow quickly from this that $\vecQ(t) = 0$ for $t\ge \gamma |\vecQ(0)|$, with $\gamma$
not depending on the particular fluid model solution; when this behavior occurs, the fluid
model is said to be \emph{stable}.

%We provide some bounds and show there is sufficient variability for our fluid process to convergence to zero. We thus show fluid model is stable. 

%As in \cite{Br08}, we say a fluid model is \emph{stable} if there exists a constant $\gamma>0$ such that  
%\begin{equation}\label{fluidstable}
%|\vecQ(t)|=0,\qquad\qquad \forall t \geq \gamma |\vecQ(0)|.
%\end{equation}

The following lemma shows that, for subcritical networks, queue sizes must vary over time and hence there is a ``mismatch" between arrival rates and departure rates at some of the queues. As we will see, this together with the relative entropy terms in Proposition \ref{HDiff} will show the Lyapunov function $H(t)$ decreases to $0$ at a uniform rate.
The lemma is proven in Subsection \ref{AppendixD} of the Appendix.  
%In Proposition \ref{pr:der.H.boud.away.0},  this bound will be used to show that the derivative is bounded away from 0.

\begin{lemma} \label{bound.change.Q}
Assume that the arrival rate vector $\veca \in \mC$.  Then
there exist constants $c,  \kappa>0$ such that, whenever $|\vecQ(0)|>0$, 
%some $j\in\mJ$ 
%and some time $t < c \, |\vecQ(0)|$
\begin{equation}
\label{eqnewkappa}
   |Q_j(t) - Q_j(0)| \geq \kappa \, |\vecQ(0)|\ 
\end{equation}
for some $j\in\mJ$ 
and $t \le c \, |\vecQ(0)|$.
\end{lemma}
\iffalse
\begin{lemma}\label{lm:log.bound}
For $x,\, y \in (0,K]$ 
\begin{equation}\label{eq:log.bound}
y\log \left(\frac{y}{x}\right) - (y-x) \geq \frac{1}{2 K} \left(y-x \right)^2.
\end{equation}
\end{lemma}
\proof{Proof}
Let 
$F(z) = z\log z - (z-1)$, with $F\rq{}(z)= \log z$ and $F\rq{}\rq{}(z) = z^{-1}$. Taking a Taylor expansion around $1$, with $F(1)=0$, for some $\theta \in(1,z)$, 
$$
F(z) = F\rq{}(1) (z-1) + \frac{1}{2} F\rq{}\rq{}(\theta) (z-1)^2 
= \frac{1}{2\theta} (z-1)^2 \geq \frac{1}{2 (1\vee z )} (z-1)^2 \ .
$$
The last equality, holds since $\theta \leq (1 \vee z)$. 
Noticing that the left hand side of \eqref{eq:log.bound} is equal to $x\, F(y/x)$ the result follos by the following relation
$$
 x F(y/x) 
 \geq \frac{x}{2 (1\vee \frac{y}{x} )} \left(\frac{y}{x}-1\right)^2
= \frac{1}{2 (x\vee y) } (y-x)^2 
 \geq \frac{1}{2 K} (y-x)^2 \ .
$$
\Halmos\endproof
\fi

\iffalse
\begin{lemma}\label{lm:H.bound}
For appropriate constants $h_1, \, h_2>0$, 
\begin{equation}\label{eq:H.bounds}
h_1 |\vecQ(t)|
\leq H(t) \leq 
h_2 |\vecQ(t)|
\end{equation}
\end{lemma}
\proof{Proof}
\Halmos\endproof
\fi

%We will use this lemma, as well as Lemma \ref{lm:log.bound} in the Appendix,  to give a lower bound on the rate of decrease of $H(t)$.

The following elementary lemma is also proved in the Appendix.

%We are now ready to prove the main result about the convergence of the Lyapunov function $H(t)$ to zero.

\begin{lemma}\label{lm:log.bound}
For $x,\, y \in (0,K]$ and any $K>0$,
\begin{equation*}
y\log \left(\frac{y}{x}\right) - (y-x) \geq \frac{1}{2 K} \left(y-x \right)^2.
\end{equation*}
\end{lemma}

We will use the two preceding lemmas to give the following lower bound on the rate of decrease of $H(t)$.  The basic argument for the proposition is similar to that leading up to (4.24) in \cite{Br96a},
although the reasoning differs in a few ways.  

\begin{proposition}\label{pr:der.H.boud.away.0}
Assume that the arrival rate vector $\veca$ belongs to the subcritical region $\mC$. 
For $T$ as in Corollary \ref{H.deriv.a.e.},  
there exist constants $c_1, c_2 >0$ such that, for $t\geq T|\vecQ(0)|$,
\begin{equation}\label{eq:H.dec.bound}
% H(u) \leq ( 1 - \epsilon ) \, H(t)\ .
  H( t + c_1\, |\vecQ(t)|) - H(t) \le - c_2  \, |\vecQ(t)|  \,  .
\end{equation}
\end{proposition}
\proof{Proof}
%The basic argument is similar to that leading up to (4.24) in \cite{Br96a},
%although the reasoning differs in a couple ways.  
 %
%The following Taylor expansion bound given in Lemma \ref{lm:log.bound} in the Appendix,
%\begin{equation}
%y\log \left(\frac{y}{x}\right) - (y-x) \geq \frac{1}{2 K} \left(y-x \right)^2
%\end{equation}
%for $x,y\in (0,K]$. We 
%

Let $K$ be the Lipschitz constant bounding the processes $D'_r(t),$ $A'_r(t),$ $D'_j(t),$ and $A'_j(t),$ over all $r\in\mR$ and $j\in\mJ$. From Proposition \ref{HDiff} and Lemma \ref{lm:log.bound}, 
\begin{align*}
H'(t) 
%= & - \sum_\rInR  D'_{l(r)r}(t) \log \frac{D'_{l(r)r}(t)}{a_r} 
%- \sum_\jInJ A'_j(t) \, \log \frac{A'_j(t)}{D'_j(t)} \\
\leq & - \sum_\rInR \frac{1}{2 K} \left(D'_r(t)-A'_r(t) \right)^2 
- \sum_\jInJ \frac{1}{2 K} \left(A'_j(t)-D'_j(t) \right)^2\\
& - \sum_\rInR (D'_r(t)-A'_r(t))  - \sum_\jInJ Q'_j(t)\\
= & - \sum_\rInR \frac{1}{2 K} \left(D'_r(t)-A'_r(t) \right)^2 
- \sum_\jInJ \frac{1}{2 K} \left(Q'_j(t)\right)^2 
\leq - \frac{1}{2 K} \sum_{j\in\mJ} \left(Q'_j(t)\right)^2  .
\end{align*}
%
%where in the equality, like in \eqref{eq:cons.flow}, we again used the conservation law of flows. Integrating is permissible by Corollary \ref{H.deriv.a.e.} and yields
%
On account of \eqref{eq:upp.Lipschitz} in Corollary \ref{H.deriv.a.e.}, this implies
\begin{equation}
\label{eqquadraticbd}
H(t) - H(s) \leq -\frac{1}{2 K}  \sum_{j\in\mJ} \int_s^t \left(Q'_j(u)\right)^2 du. 
\end{equation}

By Lemma \ref{bound.change.Q}, there exists a $j\in\mJ$ and a value of $t$ with $t \leq s + c_1 \, |\vecQ(s)|$ such that
\begin{equation*}
 |Q_j(t) - Q_j(s)| \geq \kappa \, |\vecQ(s)|.
\end{equation*}
It therefore follows from the Cauchy-Schwarz Inequality that

\begin{align}
\kappa \, |\vecQ(s)| &
%\leq |Q_j(t) - Q_j(s)| 
\leq \int_s^{s+c_1 \, |\vecQ(s)|} \left| Q'_j(u) \right| \, du \leq \left( \int_s^{s+c_1 \, |\vecQ(s)|} \left( Q'_j(u) \right)^2  du \right)^{\frac{1}{2}} \, \sqrt{c_1 \, |\vecQ(s)|}  \nonumber
\end{align} 
for this choice of $j$, and hence
$$ 
\int_s^{s+c_1 \, |\vecQ(s)|} \left( Q'_j(u) \right)^2 du \geq \frac{\kappa^2}{c_1} |\vecQ(s)|\ .
$$
Applying the last inequality to  \eqref{eqquadraticbd}, one obtains
\begin{equation*}
H(s+c_1 \, |\vecQ(s)|) - H(s) \leq - \frac{\kappa^2}{2Kc_1} |\vecQ(s)| \,,
\end{equation*}
and so   (\ref{eq:H.dec.bound}) follows by setting $c_2=\kappa^2/2Kc_1$.
\iffalse
By \eqref{eq:H.bounds} we have  
\begin{equation*}
H(t+c_1 \, |\vecQ(t)|) - H(t) \leq - \frac{K_5}{c_1}  \frac{1}{h_1} \, H(t)\ ,
\end{equation*}
and by letting $ \epsilon = K_5 / (h_1 \, c_1 )$ we have that
$H(t+c_1 \, |\vecQ(t)|) \leq ( 1 - \epsilon ) \, H(t)$.
Finally since the function $H(t)$ is non-increasing and non-negative it must be that 
\begin{equation*}
H(u) \leq ( 1 - \epsilon ) \, H(t).
\end{equation*}
for all $u \geq t + c_1 |\vecQ(t)|$.
\fi
\Halmos\endproof

Employing Lemma \ref{prop:H.not.zero}, Corollary \ref{H.deriv.a.e.}, and Proposition \ref{pr:der.H.boud.away.0}, we now show our main result for proportional switch
fluid models.  

\begin{theorem} \label{FluidStable}
If $\veca \in \mC$, then the proportional switch fluid model is stable.
\end{theorem}
\proof{Proof}
The proof follows the reasoning in \cite{Br96a}, which we repeat for completeness.
Let $c_1$ and $c_2$ be as in in Proposition \ref{pr:der.H.boud.away.0}, set $t_0=T|\vecQ(0)|$ for $T$ as in Corollary \ref{H.deriv.a.e.}, and define
\begin{equation}\label{eq:t_i.rec}
t_{i+1} = t_i + c_1 \, |\vecQ(t_i)|. \ 
\end{equation}
By Proposition \ref{pr:der.H.boud.away.0}, since $H(t)\ge 0$, 
%(Lemma \ref{prop:H.not.zero}) imply
$$
H(t_0) \geq
H(t_0) - H(t_{i+1}) \geq
\sum_{l=0}^i (H(t_l) - H(t_{l+1})) \geq \sum_{l=0}^i c_2 |\vecQ(t_l)| 
= \sum_{l=0}^i \frac{c_2}{c_1} (t_{l+1}-t_l) 
= \frac{c_2}{c_1} (t_{i+1} -t_{0})\,,
$$
and so, as $i \to \infty$, one has
$t_\infty =\lim_{i\rightarrow\infty} t_i \leq c_1 \, H(t_0)/c_2 + t_0 < \infty$.
By the continuity of $\vecQ(t)$  and the definition of $t_{\infty}$, 
$\vecQ(t_\infty)=0$, and consequently, since $H(t)$ is non-increasing,
\begin{equation*}
H(t)=0 \quad \mbox{ for } t\geq \frac{c_1}{c_2}  H(t_0) + t_0.
\end{equation*}
On the other hand, by Corollary \ref{H.deriv.a.e.}, $H(t_0)\leq \kappa_H t_0= \kappa_H T |\vecQ(0)|$.
Setting $c_3=c_1 \kappa_H T/c_2+T$, this together with Lemma \ref{prop:H.not.zero} implies 
\begin{equation}\label{eq:Q.hitting.time.of.zero}
\vecQ(t)=0 \quad \mbox{ for } t\geq c_3 \, |\vecQ(0)|\,,
\end{equation}
as required.
\Halmos\endproof

The main result in the paper, Theorem \ref{mainthrm}, follows immediately from Theorem  \ref{FluidStable} and the following proposition.
\begin{proposition}
\label{propFMQN}
Suppose that the proportional switch fluid model is stable.  Then the corresponding proportional switched network is positive recurrent.
\end{proposition} 

We postpone the proof of Proposition \ref{propFMQN} to Subsection \ref{Appendix F} of the Appendix.  We  will follow there the approach of \cite{Da95} as presented in \cite{Br08}, although the argument is considerably shorter in the present setting and requires only several pages.

\bibliographystyle{apalike}
\bibliography{references-21072016-arxiv}

\appendix

\section{Lemmas from Section \ref{LyaDeriv}}

The main aim of this section is to prove Lemma \ref{L.deriv.a.e.}. We first present the supporting Lemmas \ref{SigCont} - \ref{pr:D'.positive}.  Lemma \ref{SigCont} 
is proved in  Lemma A.3 of \cite{KeWi04}.

\begin{lemma}\label{SigCont} 
The function $\vecQ \mapsto \sigma_j(\vecQ)$ is continuous on the set $\{ \vecQ : Q_j>0\}$.
\end{lemma}

\iffalse
The $D_j(t)$ are Lipschitz, then by  \eqref{eq:ADk}
also the $A_k(t)$ are Lipschitz, whenever the classes are not source classes, 
finally for the source classes the Lipschitz property it is due to the 
i.i.d. assumption and the finite means.
\fi

\iffalse
\begin{proposition}\label{deriv.D}
\begin{subequations}
\begin{align}\label{PF.forD}
\vecD'(t) &\in \argmax_{\sigma \in \mC} \sum_{j\in\mJ} Q_j(t) \log \sigma_j & t - \mbox{a.e.} \\
D'_j(t) &= A'_j(t) & \mbox{if } Q'_j(t)=0
\end{align}
\end{subequations}
\end{proposition}
\fi

When $Q_j/|\vecQ|$ is bounded away from $0$, Lemma \ref{lm:lower.bound.sigma} states that the corresponding
coordinate $\sigma_j(\vecQ)$ of the proportional fair optimization in (\ref{PF}) is also bounded away from $0$.
%
%Lemma \ref{lm:lower.bound.sigma} bounds away from $0$ a given coordinate of the proportional fair optimization $\vecsigma(\vecQ)$ to (\ref{PF}) when the corresponding coordinate of $\vecQ$ is bounded away from $0$.  
%
%(Note that the lemma applies to general  $\vecQ$ after scaling.)
%
\begin{lemma}\label{lm:lower.bound.sigma}
Let $\vecsigma(\vecQ)$ be a solution to the proportional fair optimization (\ref{PF}), where $|\vecQ| >0$.
%
%\begin{equation}\label{pfq}
%\vecsigma(\vecq) \in \argmax_{\vecs \in \mC} \sum_{j \in\mJ } q_j \log s_j 
%\end{equation}
%where $\vecq\geq0$ and $\sum_{j} q_j = 1$. 
%
Then, for any $\epsilon>0$, there exists a $\delta>0$ such that, for any $j$, 
 $\sigma_j(\vecQ)\geq \delta$ whenever $Q_j \geq \epsilon |\vecQ|$.
\end{lemma}
\proof{Proof}
%Note the optimum of \eqref{pfq} is bounded below, for instance, by taking a constant choice of $\vecs\in\mC$. We proceed by contradiction.
Scaling $\vecQ$ by $|\vecQ|$, it suffices to prove the result for $\sum_{j'\in \mJ}Q_{j'} = 1$, when $Q_j \ge \epsilon$.  If the result were not true, then, for some $\epsilon>0$ and $j$, there would exist a sequence $\vecQ^{(k)}$ with $Q_j^{(k)}>\epsilon$ and  $\sigma_j(\vecQ^{(k)}) \rightarrow 0$ as $k\rightarrow\infty$.
%; by taking an appropriate subsequence, one can assume $\vecQ^{(k)}$ also converges. 
It would follow that
\begin{equation*}
\sum_{j' \in\mJ } Q^{(k)}_{j'} \log \sigma_{j'}(\vecQ^{(k)}) \leq Q_j^{(k)} \log \sigma_j(\vecQ^{(k)}) + (1- Q_j^{(k)}) \log \sigma_{\text{max}}
%
%\sigma_j(\vecQ^{(k)}) + (1- Q_j^{(k)}) \log \max_{j' \in\mJ } \{\sigma_{j'}\} 
%
 \xrightarrow[k\rightarrow\infty]{} -\infty,
\end{equation*}
where $\sigma_{\text{max}}$ is defined at the beginning of Subsection \ref{subsection2.1}.
Since the maximum in (\ref{PF}) is bounded from below (by any fixed choice of $\vecsigma$), this gives a contradiction.
%
\iffalse
Consider the value
$$ s(\epsilon) = \inf_{\vecalpha:q_{j'} \geq \epsilon} \hat\sigma_{j'} $$
and assume that it is equal to $0$. The we can find a sequence $\vecalpha^{[k]}$
such that the corresponding optimal solutions $\hat\vecsigma^{[k]}$ have the $j'$ component vanishing.

However we would have that 
\begin{align}
q^{[k]}_{j'} \log \sigma^{[k]}_{j'} + \sum_{\tiny\begin{array}{c}j \in\mJ \\ j\not=j'\end{array}} q^{[k]}_j \log \sigma^{[k]}_j 
\leq
\epsilon \log \sigma^{[k]}_{j'} + |\mJ| \sigmamax \to -\infty
\end{align}
that is a contradiction, since $\hat\vecsigma^{[1]}$ would be a feasible solution with higher value of the objective function
$$\inf_{\vecalpha:q_{j'} \geq \epsilon} \left(
 q_{j'} \log \hat\sigma^{[1]}_{j'} 
+ \sum_{\tiny\begin{array}{c}j \in\mJ \\ j\not=j'\end{array}} q_j \log \hat\sigma^{[1]}_j 
\right) 
\geq 
|\mJ| \log \left( \min_{j \in\mJ}  \hat\sigma^{[1]}_j \right).$$
\fi
\Halmos\endproof

%This lemma proves once all queues have emptied there initial work they will have non-zero departure rate. In comparison to \cite[Lemma 5]{Ma07}, this lemma helps us characterize the derivative of $L(t)$ with equality rather than an inequality bound.
%
Lemma \ref{pr:D'.positive} bounds $D_k(t)/t$ from below for large $t$.  A related result is given in Lemma 5 of \cite{Ma07}.

\begin{lemma}\label{pr:D'.positive} %{pr:D.unbounded}
There exists $T >0$, not depending on $\vecQ(0)$,  such that, for any class $k$ and $t\ge T|\vecQ(0)|$,  
\noindent i) for appropriate $\beta >0$,
 $D_k(t) >   \beta t$  and \noindent ii)  $D'_k(t)>0$ a.e.

%\begin{lemma}\label{pr:D'.positive} %{pr:D.unbounded}
%\noindent i) There exists $\beta>0$ such that  
% $$D_k(t) >   \beta t, \qquad t\geq 0.$$ 

\iffalse
For each $k\in\mK$ 
  \begin{equation}\label{eq:D.unbounded}
  \limsup_{t\rightarrow\infty} \frac{D_k(t)}{t} > 0\, .
  \end{equation}
\fi
%\noindent ii)
%There exists $T_j>0$ such that
%  \begin{equation}\label{eq:D'.positive}
%  D'_j(t)>0, \qquad a.e. \quad t \geq |\vecQ(0)| T_j  \ . 
%  \end{equation}
\end{lemma}
\proof{Proof}
%Both proofs involve reasonably straightforward induction arguments. 
%We give a brief proof of both results.

\noindent i) We argue by induction along each route, assuming for a class $k$, times $s$ with $s\ge T_1|\vecQ(0)|$ for a given $T_1\ge 1$, and a given $\alpha>0$, that $A_k(s)\geq \alpha s$.
We first show the analog of the desired result, but for $D_j(t)$ instead of $D_k(t)$, where $k\in j$.
If $D_j(s) < \alpha s/2$, then $Q_j(s)\geq Q_k(s)\geq \alpha s/2 $, and so, when $s\ge T_1|\vecQ(0)|$,
%
%Thus, as work arrives linearly into the network, a positive proportion of the work in the network is at queue $j$, 
%
\begin{equation*}
\frac{Q_j(s)}{\sum_{j\rq{}\in\mJ} Q_{j\rq{}}} \geq \frac{\alpha s/2 }{\sum_{j'\in\mJ} (a_{j'} s + Q_{j'}(0)) } \geq \epsilon
\end{equation*}
for appropriate $\epsilon>0$. 
Therefore, by Lemma \ref{lm:lower.bound.sigma}, $D'_j(s) \geq \delta$, where $\delta$ does not depend on the initial queue size distribution; from this, it follows
%
%So either $D_j(t)$ lies above the line $\alpha/2 \cdot t$ or $D_j(t)$ moves above the linear rate $\delta$. 
%
that, if $D_j(s) < \alpha s/2$ for all $s\in [t/2,t]$ with $t\ge 2T_1|\vecQ(0)|$, then $D_j(t)\ge \delta t/2$.
Combining the last inequality with the case when $D_j(s) \ge \alpha s/2$ for some $s\in [t/2,t]$,  and setting $ \delta'= \min\{ \alpha/4 , \delta/2 \}$, one obtains $D_j(t) \geq  \delta't$ for $t\ge 2T_1|\vecQ(0)|$ in both cases.

We now show $D_k(t)$ also increases linearly in time.  Note that $A_j(t) \le \bar{a} t $, where $\bar{a}= \max_j{a_j}+ |\mK| |\sigma_{\max}|$. Applying this and setting $\gamma=  \delta'/2\bar{a}$, one can check that 
$D_j(t) \geq A_j(\gamma t) + |\vecQ(0)|$ for $t\ge T_2|\vecQ(0)|$, with $T_2 = 2 \max(T_1, 1/\delta')$ .
Using \eqref{eq:Dk} and the FIFO assumption \eqref{eq:Ak}, one therefore obtains
\begin{equation*}
D_k(t) = \Gamma_k(D_j(t)) \geq \Gamma_k(A_j(\gamma t) + Q_j(0)) = A_k(\gamma t) + Q_k(0) \geq \gamma \alpha t
\end{equation*}
for $t\ge T_2|\vecQ(0)|$, where the last inequality follows from the induction assumption.  Since $A_k(t)=a_r t$ for the first class on a route $r$, applying (\ref{eq:ADk}) and the argument just given, we can inductively show that for all classes on a given route $r$, $D_k(t) > \beta t$ for $t \ge T|\vecQ(0)|$, for appropriate choices of $\beta,T>0$. This is the desired result.

\iffalse
Seeking a contradiction, suppose $\lim_{t\rightarrow\infty} D_k(t)/t = 0$ and $\limsup_{t\rightarrow \infty} A_k(t)/t >0$ for some class $k$ at some queue $j$. %There is a positive arrival rate $a_k$ to that class. So since arrival rates to queue $j$ are bounded above, there must be a positive density of class $k$ work in queue $j$. 
Then $\limsup_{k\rightarrow\infty} Q_k(t)/t >0$, thus a positive proportion of the work in the network must be class $k$ work distributed along queue $j$. $Q_j(t)$ is greater than $Q_k(t)$. By the previous lemma and since $Q_j(t)$ is Lipschitz, $D\rq{}_j(t)$ must be bounded below for a positive proportion of time. Thus $\limsup_{t\rightarrow\infty} D_j(t)/t>0$. Further since service is FIFO and arrival rates to queue $j$ are bounded, over time positive proportions of class $k$ arrivals must depart. So one sees that $\limsup_{t\rightarrow\infty} D_k(t)/t>0$. This yields a contradiction. 

Clearly $\limsup_{t\rightarrow \infty} A_k(t)/t>0$ holds for each input class.  
The remainder of the proof applies the argument above by induction along the classes of each route.
\fi
\noindent ii) We again argue by induction along each route, assuming for a class $k$ that
%%We consider the times $t$ were the processes arrival and departure processes are differentiable, which holds almost everywhere.
%Suppose that, for class $k$ at queue $j$, 
%
$A'_k(t) > 0$ a.e. on $t\ge T_3|\vecQ(0)|$ for some $T_3\ge T$, where $T$ is given in part i); we henceforth restrict our attention to times $t$ where the arrival and departure processes are differentiable.  It follows under this assumption that 
\begin{equation}\label{Ddash}
D'_j(t)>0 \qquad \text{for } t\ge T_3|\vecQ(0)|
\end{equation}
because, if $Q_j(t) > 0$, then $D_j'(t)>0$ by Lemma \ref{lm:lower.bound.sigma}, whereas, if $Q_j(t)=0$, then $Q_j'(t)=0$ and so $D_j'(t)=A_j'(t)\geq A_k'(t)>0$. 

By \eqref{eq:Ak} and the definition of $\Gamma'_k$,  for $t\ge T_3|\vecQ(0)|$,
\begin{equation}
\label{eqnew42}
\Gamma'_k (A_j(t)+ Q_j(0)) = \frac{A'_k(t)}{A'_j(t)} > 0\,.
\end{equation} 
Also, by part i), for  $t \ge T_4 |\vecQ(0)|$, with $T_4 = \max(T , \beta^{-1} (A_j(T_3 |\vecQ(0)|)/|\vecQ(0)| +1))$, 
%
%over $A_j(T_0 |\vecQ(0)|)$ units of work have been served at queue $j$, that is 
%
one has $D_j(t) \geq A_j(T_3 |\vecQ(0)|) + Q_j(0)$; together with (\ref{eqnew42}), this implies that 
\begin{equation}\label{Gdash}
\Gamma'_k(D_j(t))>0 
\end{equation}
for $t\geq T_4 |\vecQ(0)|$.
Using the same upper bound for $A_j(t)$ as in part i), one can show that $T_4 \le T_5$ for $T_5 = \max(T_3 , \beta^{-1}(\bar{a}T_3 + 1))$.
So, \eqref{Ddash} and \eqref{Gdash} together imply that 
\[
D_k'(t)= D_j'(t) \Gamma'_k(D_j(t))>0
\] 
for $t \ge T_5 |\vecQ(0)|$. Applying (\ref{eq:ADk}), one can then inductively repeat this argument along the classes of each route to give the desired result after a new choice of $T$.
 \Halmos
\endproof

We now employ Lemmas \ref{SigCont} and  \ref{pr:D'.positive} to demonstrate Lemma \ref{L.deriv.a.e.}.

\proof{Proof of Lemma \ref{L.deriv.a.e.}} 
i)  We first obtain lower and upper bounds for  $L(t+h) - L(t)$.
The processes $\vecD(t)$ and $\vecQ(t)$ are almost everywhere differentiable and,  by part ii) of Lemma \ref{pr:D'.positive}, there exists a $T>0$ such that $D_j'(t) >0$ a.e. on $t \geq T|\vecQ(0)|$ for each queue $j$; choose a $t$ satisfying these conditions.
\iffalse
From \eqref{def:entropy}, $H(t)$ admits derivative as soon as the term 
\begin{equation}
L(t) = \sum_{j\in\mJ} Q_j(t) \log D'_j(t)
\end{equation}
does.
\fi
%
Since $\vecD'(t)$ is suboptimal for the proportional fair optimization (\ref{PF}) with queue size vector $\vecQ(t+h)$ for given $h$,
\begin{equation}
\label{eqnew44}
\begin{split}
L(t+h) - L(t) &= \sum_\jInJ \left(Q_j(t+h)\log D'_j(t+h)-Q_j(t)  \log D'_j(t)\right) \\
&\geq \sum_\jInJ \left(Q_j(t+h)-Q_j(t)\right)  \log D'_j(t).  
\end{split}
\end{equation}
Also, since $\vecD'(t+h)$ is suboptimal for queue size vector $\vecQ(t)$, 
\begin{equation}
\label{eqnew45}
\begin{split}
L(t+h) - L(t) 
\leq  & \sum_{j\in\mJ} \left(Q_j(t+h)-Q_j(t)\right)  \log D'_j(t+h) \\
\leq  & \sum_{\substack{j\in\mJ:\\ Q_j(t) >0}} \left(Q_j(t+h)-Q_j(t)\right)  
\log D'_j(t+h) +\sum_{\substack{j\in\mJ:\\ Q_j(t) =0}} Q_j(t+h)
\log \sigma_{\max} .
\end{split}
\end{equation}
(After splitting the sum over $j$ into two parts depending on whether or not $Q_j(t)=0$, we employed $\log D_j'(t) \leq \log \sigma_{\max}$ in the summation over $j$ with $Q_j(t)=0$.)

After dividing by $h$, we consider the limiting behavior of the left hand sides of (\ref{eqnew44}) and (\ref{eqnew45}), as $h\rightarrow 0$, by employing the bounds on the right hand sides of the equations.
Since $\vecQ'(t)$ exists, after dividing by $h$, the right hand side of  (\ref{eqnew44}) converges to
\begin{equation}\label{Diffy}
\sum_\jInJ Q'_j(t) \log D'_j(t) 
\end{equation}
as  $h\rightarrow 0$. On the other hand,  by Lemma \ref{SigCont},  $D'_j(t)=\sigma_j(\vecQ(t))$ is continuous where $Q_j(t)>0$, and so $D'_j(t+h) \rightarrow D'_j(t)$ as $h\rightarrow 0$. Recalling that $Q'_j(t)=0$ when $Q_j(t)=0$, the first term on the right hand side of (\ref{eqnew45}) also converges to (\ref{Diffy}) after dividing by $h$. 
Moreover, the last term in (\ref{eqnew45}) converges to $0$ after dividing by $h$, since
$Q'_j(t)=0$ for such $j$.  Putting these limits together, (\ref{eq:EnvelopeTheorem}) follows.
%\begin{equation*}
%L'(t) = 
%\sum_\jInJ Q'_j(t) \log D'_j(t),
%\end{equation*}
%as required.

ii) Most of the work consists of computing an upper bound for the upper right Dini derivative of the function $L(\cdot)$  at time $t$, which is given by
\begin{equation}\label{eq:L.upper.Dini.deriv}
D^+L(t) = \limsup_{h\searrow 0} \frac{ L(t+h) - L(t)}{h} \, .
\end{equation}
By the same reasoning as in (\ref{eqnew45}), for
$h>0$,
\begin{equation}
\label{eqnew48'}
\begin{split}
\frac{L(t+h) - L(t)}{h} 
&\leq  \sum_{j\in\mJ} \frac{Q_j(t+h)-Q_j(t)}{h} \log D'_j(t+h) \\
& \leq  
\sum_{\substack{j\in\mJ:\\ Q_j(t) =0}} \frac{Q_j(t+h)}{h} \log D'_j(t+h)  
 + \sum_{\substack{j\in\mJ:\\ Q_j(t) >0}}  \frac{Q_j(t+h)-Q_j(t)}{h} \log D'_j(t+h).
\end{split}
\end{equation}
Denoting by $K_A$ and $K_Q$ the Lipschitz constants for the processes $\vecA(t)$ and $\vecQ(t)$, and employing $\vecQ(t) = \vecA(t)-\vecD(t)$, it follows that this is at most
\begin{equation*} 
(K_Q +K_A) |\mJ| \log \sigmamax
- \sum_{\substack{j\in\mJ:\\ Q_j(t) >0}}  \frac{D_j(t+h)-D_j(t)}{h} \log D'_j(t+h). 
%
%\sum_{\substack{j\in\mJ:\\ Q_j(t) =0}} K_Q \log \sigmamax
%+ \sum_{\substack{j\in\mJ:\\ Q_j(t) >0}}  K_A \, \log \sigmamax
%- \sum_{\substack{j\in\mJ:\\ Q_j(t) >0}}  \frac{D_j(t+h)-D_j(t)}{h} \log D'_j(t+h). 
%
\end{equation*}
By Lemma \ref{SigCont} , $D'_j(t+h)$ is continuous at $h=0$ when $Q_j(t)>0$.  Taking the $\limsup$ of the left hand side of (\ref{eqnew48'}), it follows from the above inequalities that
\begin{align}\label{eq:upper.right.derivative}
D^+L(t) \leq  
(K_Q +K_A) |\mJ| \log \sigmamax  - \sum_{\substack{j\in\mJ:\\ Q_j(t) >0}}  D'_j(t) \log D'_j(t) \leq \,  \kappa_L \ 
%
%\sum_{\substack{j\in\mJ:\\ Q_j(t) =0}} K_Q \log \sigmamax
%+ \sum_{\substack{j\in\mJ:\\ Q_j(t) >0}}  K_A \, \log \sigmamax \nonumber  - \sum_{\substack{j\in\mJ:\\ Q_j(t) >0}}  D'_j(t) \log D'_j(t) \leq \,  \kappa_L \ 
%
\end{align}
for some constant $\kappa_L$, since  $x\,\log x$ is bounded from below. It then follows from standard results on Dini derivatives that
\begin{equation}\label{LDev}
L(t)-L(s) \leq \kappa_L (t-s)
\end{equation}
%
% see Theorem 2.3, page 348 in [[CITE BOOK ON DINI DERIVATIVE]] % BERNARDO: I was not able to find the author and the title of this book
(see, for instance, Theorem 3.4.5 on page 65 of \cite{kannan1996advanced}), which is the desired inequality.

iii) For $n\in\bN$ and fixed $s$, $t$, and $u$, with $u\in [s,t]$, define the dyadic floor and ceiling  $\lfloor u \rfloor_n$ and $\lceil u \rceil_n$ of $u$ by 
\begin{align*}
\lfloor u \rfloor_n &:= \max_{k\in\bZ_+}\{ k (t-s) 2^{-n} + s :  k (t-s) 2^{-n} + s \leq u\}, \\
 \lceil u \rceil_n &:=\min_{k\in\bZ_+}\{ k (t-s) 2^{-n} + s :  k (t-s) 2^{-n} + s > u\};
\end{align*}
note that $\lceil u \rceil_n -\lfloor u \rfloor_n= (t-s)2^{-n}$. By interpolating terms and rewriting them in integral form, one has
\begin{align*}
 L(t) - L(s)
 = & \sum_{k=1}^{2^{n}} \Big( L (k(t-s)2^{-n}+s)- L ((k-1)(t-s)2^{-n}+s) \Big) \\
 =  & \int_s^{t} \frac{ L (\lceil u \rceil_n)- L (\lfloor u \rfloor_n)}{\lceil u \rceil_n -\lfloor u \rfloor_n}  \, du 
  =   \lim_{n \rightarrow \infty} \int_s^t\frac{ L (\lceil u \rceil_n)- L (\lfloor u \rfloor_n)}{\lceil u \rceil_n -\lfloor u \rfloor_n}  \,  du \\
\leq & \int_s^t \limsup_{n \rightarrow \infty}   \frac{ L (\lceil u \rceil_n)- L (\lfloor u \rfloor_n)}{\lceil u \rceil_n -\lfloor u \rfloor_n}  \,  du 
=   \int_s^t L'(u) du.
\end{align*}
The limit in the third equality holds trivially since the expression is not a function of $n$; the inequality follows from the Fatou Inequality since the integrand is bounded above by (\ref{eq:L.Lipschitz.up}); and the last equality follows by using (\ref{eq:EnvelopeTheorem}).
\Halmos\endproof

\iffalse
Defining 
 $\bar f_n(x) = f_n(\left\lceil 2^n \, (x-s)/(t-s) \right\rceil)$ 
and
 $\underline f_n(x) = f_n(\left\lfloor  2^n \, (x-s)/(t-s)  \right\rfloor)$ 
it follows that 
\begin{align}
 L(t) - L(s)
 =   \int_s^t \frac{L (\bar f_n(x))- L (\underline f_n(x))}{(t-s)2^{-n}}   \,  dx 
 = &  \lim_{n \rightarrow \infty} \int_s^t \frac{L (\bar f_n(x))- L (\underline f_n(x))}{(t-s)2^{-n}}   \,  dx \nonumber 
\end{align}
\fi
\iffalse
Since $| \bar f_n(x) - \underline f_n(x)| \leq 2^{-n}$, using the upperbound \eqref{eq:L.Lipschitz.up} and the bounded convergence theorem we get
\begin{align}
 L(t) - L(s)
\leq  &  \int_s^t 
 \lim\sup_{n \rightarrow \infty} 
\frac{L (\bar f_n(x))- L (\underline f_n(x))}{\bar f_n(x) - \underline f_n(x)} \, 
\frac{\bar f_n(x) - \underline f_n(x)}{(t-s)2^{-n}}   \, dx  \nonumber \\
= & \int_s^t \frac{dL(x)}{dt} dx.
\end{align}
where we used that $\bar f_n(x), \, \underline f_n(x) \to x$ and the second fraction converges to $1$.
\fi

\section{Stability Results from Section \ref{PosRec}} \label{AppendixD}

For the proof of Lemma \ref{bound.change.Q}, we require the following technical lemma. 
%[ CAN THIS LEMMA BE MERGED WITH LEMMA \ref{lm:lower.bound.sigma}  ]

\begin{lemma}\label{lm:lowerbound.q.sigma}
For any $\veca\in\mC$, there exist $\epsilon >0$ and $c>0$, such that, for any $\vecQ\in\bR_+^{|\mJ|}$ with $|\vecQ|>0$, 
 one has $Q_j \geq c \, |\vecQ|$ and $\sigma_j(\vecQ) \ge (1+\epsilon) a_j $ for some queue $j$.
\end{lemma}
\proof{Proof}
Since the assertion is trivial if $a_j =0$ for some $j$, we can assume that $a_j >0$ for all $j$.

Since $\veca \in \mC$, there exists an $\epsilon>0$ such that 
$(1+\epsilon)^2 \, \veca \in  \mC$. Moreover, since $(1+\epsilon)^2 \veca$ is  suboptimal, 
\begin{equation*}
\sum_\jInJ Q_j \log\sigma_j(\vecQ) 
\geq \sum_\jInJ Q_j \log \left( (1+\epsilon)^2 \, a_j \right);
%
%&= \sum_\jInJ Q_j \log \left( (1+\epsilon) \, a_j \right)
%  + \sum_\jInJ Q_j \log (1+\epsilon) \ .
%
\end{equation*}
hence
$$
\sum_\jInJ Q_j \log \frac{\sigma_j(\vecQ)}{(1+\epsilon) a_j} 
\geq |\vecQ|\log (1+\epsilon),
%=  \sum_\jInJ \frac{|\vecq|}{|\mJ|} \log (1+\epsilon)
$$
which implies the existence of at least one queue, say $j$, with
$$
Q_j \log \frac{\sigma_j(\vecQ)}{(1+\epsilon) a_j} 
\geq \frac{|\vecQ|}{|\mJ|} \log (1+\epsilon) 
> 0.
$$
This implies that $\sigma_j(\vecQ) > (1+\epsilon) a_j$, as the logarithm on the right hand side is positive, and also that
$$ \frac{Q_j}{|\vecQ|} \geq \frac{1}{|\mJ|}  \log (1+\epsilon) \left(\log \frac{\sigma_j(\vecQ)}{(1+\epsilon) a_j}\right)^{-1} \geq \frac{1}{|\mJ|} \log (1+\epsilon) \left(\log \frac{\sigma_{\max}}{(1+\epsilon) a_{\min}}\right)^{-1} \ , $$
where $a_{\min} := \min \{a_j: \jInJ\}$.
Denoting by $c>0$ the right hand side of this expression, the desired lower bound on $Q_j$ follows.
\Halmos\endproof

\iffalse
\begin{lemma} \label{bound.change.Q}
There exist constants $c, \, K>0$, such that, for $|\vecQ(0)|>0$, some $j\in\mJ$ 
and some time $t < c \, |\vecQ(0)|$
\begin{equation}
   |Q_j(t) - Q_j(0)| \geq K \, |\vecQ(0)|\ .
\end{equation}
\end{lemma}
\fi

We now prove Lemma \ref{bound.change.Q}, which gives a lower bound on the amount $Q_j(t)$ will change over time for some $j$. The main idea is that, if
$Q_j(t)$ remains nearly constant, then, using
Lemma \ref{lm:lowerbound.q.sigma}, it will follow that one of the queues $j$ will empty at a faster rate than mass enters its route, which will provide a contradiction over a long enough time interval. 

%
%proof is that the proportional scheduler $\vecsigma(\vecQ)$ must be bounded above one component of the vector arrive rates $\veca=(a_j:j\in\mJ)$, since $\veca$ non-optimal for the proportional fair optimization and $\vecsigma$ is optimal. Now either that component of the schedule can dominates the arrival rate for a long period of time, in which case queue sizes associate with the queue are drained, or the value of $\vecsigma(\vecQ)$ changes, which can only occur due to a change in queue size. In either case there is a change in queue sizes over time.
%

\proof{Proof of Lemma \ref{bound.change.Q}}
For $|\vecQ(0)|>0$, it follows from  Lemma \ref{lm:lowerbound.q.sigma} that, for some queue $j$, 
$\sigma_{j}(\vecQ(0)) \ge  ( 1+ \epsilon) a_{j}$, where $\veca$ is the arrival rate vector.
Denote by $\tilde Q_{j}(0)$ the quantity of packets at queue $j$ at time $0$, together with the packets already in the network then
that will eventually be routed through $j$.
Setting $\tilde T = 2 \tilde Q_{j}(0) /(\epsilon \, a_{j})$ and $T = 2 |\vecQ(0)|/(\epsilon \, \amin)$ (with $\amin := \min_{j' \in \mJ} a_{j'}$), one has $T \geq \Tilde T$; note that $\amin > 0$. 
%
%\, \amin)$, one has $T \geq \Tilde T$. 
%

It is not possible that
\begin{equation}
\label{eqnewlemma3}
\sigma_{j}(\vecQ(t)) > a_{j} \left( 1+ \frac{\epsilon}{2} \right) \qquad \text{for all } t \in [0, T]. 
\end{equation}
For, if (\ref{eqnewlemma3}) were to hold, then
$\int_0^T  \sigma_{j}(\vecQ(t)) \, dt > a_{j} \left( 1+ \frac{\epsilon}{2} \right) \, T $,
from which the contradiction
\begin{align}
0 & \,\leq \, \tilde Q_{j}(T) 
  \, = \, \tilde Q_{j}(0) + a_{j} \, T - \int_0^T \sigma_{j}(\vecQ(t)) \, dt \nonumber  \,<\, \tilde Q_{j}(0) - a_{j}  \frac{\epsilon}{2}  \, T 
  \, \leq \, \tilde Q_{j}(0) - a_{j}  \frac{\epsilon}{2}  \, \tilde T \,=\, 0 \nonumber
\end{align}
would follow. Hence, for some $t \in [0, T]$,
$\sigma_{j}(\vecQ(t)) \leq a_{j} \left( 1+ \frac{\epsilon}{2} \right) $;
since the index $j$ was obtained from Lemma \ref{lm:lowerbound.q.sigma}, 
this implies that 
\begin{equation}\label{sigma.distance}
\sigma_{j}(\vecQ(0)) - \sigma_{j}(\vecQ(t)) \ge  a_{j} \frac{\epsilon}{2} \ . 
\end{equation}
By Lemma \ref{SigCont}, there exists $\delta > 0$ such that,
for any $\vecQ$ with $|\vecQ - \vecQ(0)| \leq \delta \, |\vecQ(0)|$, 
\begin{equation}
\label{eqafter53}
|\sigma_{j}(\vecQ) - \sigma_{j}(\vecQ(0))| \leq a_{j} \frac{\epsilon}{4}\,.
  \end{equation}
It therefore follows from \eqref{sigma.distance} that,
with possibly a new choice of $j$,
\begin{equation} \label{eq:delta}
| Q_j(t) - Q_j(0) | > \frac{\delta}{|\mJ|} \, |\vecQ(0)|\,.
\end{equation}
 
We still need to show that the bound in (\ref{eqnewkappa}) is uniform in
$|\vecQ(0)| >0$.  The choice of $\epsilon$, which is determined by Lemma
\ref{lm:lowerbound.q.sigma}, does not depend on $\vecQ(0)$.
Employing Lemma \ref{SigCont}
and  the compactness of $\{\vecQ\in\bR_+^{|\mJ|}:|\vecQ|=1\}$, $\delta$ in \eqref{eq:delta} can also be chosen so as not to depend on
$\vecQ(0)$ since $Q_j (0)/\vecQ(0) \ge c >0$ by
Lemma \ref{lm:lowerbound.q.sigma}, for $j$ in (\ref{eqafter53}).
Setting $\kappa = \delta / |\mJ|$, we obtain (\ref{eqnewkappa}).
\Halmos\endproof

The elementary Lemma \ref{lm:log.bound} follows from a straightforward Taylor expansion.

\proof{Proof of Lemma \ref{lm:log.bound}}
Setting
$F(z) = z\log z - (z-1)$, one has $F(1) = F\rq{}(1)= 0$ and $F\rq{}\rq{}(z) = z^{-1}$.  Expanding around $1$ gives 
$$
F(z) = \frac{1}{2} F\rq{}\rq{}(\theta) (z-1)^2 
= \frac{1}{2\theta} (z-1)^2 \geq \frac{1}{2 \max(1,z)} (z-1)^2 
$$
 for some $\theta$ between $1$ and $z$.
%The last equality, holds since $\theta \leq \max(1, z)$. 
Since
\[
x\, F(y/x)=y\log \left(\frac{y}{x}\right) - (y-x)
\]
for $x,y >0$, the desired inequality follows from
$$
 x F(y/x) 
 \geq \frac{x}{2 \max(1,y/x )} \left(\frac{y}{x}-1\right)^2
= \frac{1}{2 \max(x,y) } (y-x)^2 
 \geq \frac{1}{2 K} (y-x)^2 
$$
for $K\ge \max(x,y)$.
\Halmos\endproof

\section{Fluid Limit}
\label{appendFMFL}

\iffalse
We define the processes $A_k(t)$ also for $t\leq0$ such that $Q_k(0)=A_k(0)$.
For the discrete and the continuous model the following equations hold
\begin{align}
D_j(t) &= \sum_{s=0}^t \sigma_j(\vecQ(s)) \label{eq:sigmaj}\\
D_j(t) &= \int_0^t \sigma_j(\vecQ(s)) \, ds \label{eq:cont.sigmaj} \ .
\end{align}
For the discrete time process, for $t > 0$
\begin{equation}\label{eq.arr.disc.time}
A_r(t+1)-A_r(t) \quad  r\in\mR
\end{equation}
are i.i.d. random variables with mean $a_r$, while for  the continuous fluid model and $t>0$
\begin{equation}\label{eq.fluid.model}
A_r(t) - A_r(0) = a_r \, t \quad r\in\mR \ .
\end{equation}
Let $\invA_j(\cdot)$ be the inverse function of $A_j(\cdot)$ that is
\begin{equation}
\invA_j(A_j(t))=\inf\{s:A(s)\geq t\}
\end{equation}
such that $A_j(\invA_j(s))=s$.
Notice, in the second equality that we do not allow for impulsive arrivals.
We have that
\begin{equation}
\Gamma_k(s) = A_k(\invA_j(s))
\end{equation}
For the discrete time model
\begin{equation}
\gamma_k(s) = \Gamma_k(s+1) - \Gamma_k(s) 
\end{equation}
and in the continuos time model
\begin{equation}
\gamma_k(s) = \Gamma'_k(s)  
\end{equation}
\fi

In Section \ref{PROOF}, we employed Proposition \ref{propFMQN} to conclude that
the positive recurrence of a proportional switched network follows from the stability
of the corresponding proportional switch fluid model.  
One of the main steps in showing Proposition \ref{propFMQN} 
is Proposition \ref{FluidLimit}, which states that
the proportional switch fluid model in Section \ref{fluidmodel} is the limit of scaled discrete time switched networks.
We demonstrate  Proposition \ref{FluidLimit} in this subsection; our approach is similar to the standard fluid limit approaches in \cite{Da95} and \cite{Br08}.

 We consider a family of proportional switched networks indexed by $c\in\bZ_+$, employing
analogous notation 
to that introduced in Section \ref{Model}, e.g., $A^c_x$, $D^c_x$, $Q^c_x$ for $x\in \mJ \cup \mK \cup \mR$.  We will assume that the sum of the initial queue sizes of the queueing network is equal
to $c$, that is,
%\begin{equation*}
$c=|\vecQ^c(0)| = \sum_{j\in\mJ} Q^c_j(0)$.
%\end{equation*}
%%%
%%%
To simplify the proof, we couple the processes on the same probability space, assuming that the external arrival processes are identical for different $c$, that is, for $c\in\bZ_+$ and $r\in\mR$,
\begin{equation*}
 A^c_r(t) =  A_r(t).
\end{equation*}
 
For $x\in \mJ \cup \mK \cup \mR$ and $k\in\mK$, we introduce the 
scaled processes $\bar{A}^c_x(t)$, $\bar{D}^c_x(t)$, $\bar{Q}^c_x(t)$, and $\bar{\Gamma}^c_x(t)$
on $t\ge 0$, by setting
\begin{equation}\label{eq:fluid.scaled.processes}
\bar{A}^c_x(t) = \frac{A_x^c( ct )}{c},\;\;\bar{D}^c_x(t) = \frac{D_x^c(ct)}{c},\;\;\bar{Q}^c_x(t) = \frac{Q_x^c(ct)}{c},\;\; \bar{\Gamma}^c_k(t) = \frac{\Gamma^c_k(ct)}{c}  
\end{equation}
for $t\in\{0,c^{-1},2c^{-1}, 3c^{-1},...\}$, and interpolating linearly.
The following limit result holds for these processes. 
%in (\ref{eq:fluid.scaled.processes}). 
Here, $G = G_1 \cap G_2$ is a set with $\mathbb{P}(G) = 1$, where $G_1$ and $G_2$ will be specified shortly.
\begin{proposition}\label{FluidLimit}
For each $\omega \in G$, any scaled subsequence  $ (\bar{A}^{c_i}_x, \bar{D}^{c_i}_x,\bar{Q}^{c_i}_x,\bar{\Gamma}_k^{c_i} : x\in \mJ \cup\mK \cup \mR, \,k\in\mK,\,c_1 < c_2 < \ldots)$ of a sequence of proportional switched networks contains a further subsequence that converges uniformly on compact time intervals. Moreover, any such limit satisfies the proportional switch fluid model equations (\ref{eq:Incr}-\ref{eq:ADk}) and (\ref{Fluid:1}-\ref{PF}). 
\end{proposition}

In order to demonstrate Proposition \ref{FluidLimit}, we recall that, for each $r\in \mR$,  $(A_r(t)-A_r(t-1): t\in \mathbb{Z}_+)$ is a collection of i.i.d. random variables with  
mean $a_r < \infty$. It thus follows from
the Strong Law of Large Numbers that, on a set $G_1$ with $\mathbb{P}(G_1) = 1$, 
one has that, for each $r\in\mR$,  
$\bar{A}^c_r(t)\to a_r t$ as $c\to \infty$
uniformly on compact time intervals.

We also introduce, for each $j$, the martingales
\[
M_j^c(t) = \sum_{\tau=1}^t \left[ \pi_j({\vecQ}^{c}(\tau)) - \sigma_j({\vecQ}^{c}(\tau)) \right] ,
\]
 where $\pi_j(\cdot)$ and $\sigma_j(\cdot)$  are as in (\ref{DServe}), and set
$\bar{M}^c_j(t) = M^c_j(ct)/c$.
Since the increments of $\bar{M}^c_j(t)$ are uniformly bounded, it will follow from standard 
martingale bounds that, on a set with probability $1$, $\bar{M}_j^c(t) \to 0$ as $c\to \infty$, uniformly on compact time intervals for each $j\in\mJ$:

\begin{lemma} 
\label{lastlemma}
On a set $G_2$ with $\mathbb{P}(G_2) = 1$,
\[
\sup_{t\leq T} \Big|\bar{M}_j^c(t)  \Big|\xrightarrow[c\rightarrow \infty]{} 0
\]
 for all $j\in\mJ$ and $T>0$.
\end{lemma}

We will demonstrate Lemma \ref{lastlemma} after the proof of Proposition \ref{FluidLimit}. 

The proof of Proposition \ref{FluidLimit} is relatively straightforward. Most of the work consists
of showing the above limits satisfy the equations (\ref{eq:Incr}-\ref{eq:ADk}) and (\ref{Fluid:1}-\ref{PF}); the invariance of
 the equations (\ref{eq:Incr}-\ref{eq:ADk}) under fluid scaling is an important ingredient. 

%The proof is somewhat standard and can be found in the appendix of this paper. %The main interesting point is in assuring that the condition (\ref{Fluid:3}) holds.
\proof{Proof of Proposition \ref{FluidLimit}.}
%Since, for each $r\in \mR$,  $(A_r(t)-A_r(t-1): t\in \mathbb{Z}_+)$ are i.i.d. random variables with  mean $a_r < \infty$, it follows from
%the Strong Law of Large Numbers that $\lim_{t\to\infty}\bar{A}^c_r(t) = a_r$
%on a set $G$ with $\mathbb{P}(G) = 1$ (over the index set $t\in \mathbb{R}_+$).

By the Arzel\`a-Ascoli Theorem (cf. \cite{Bi99}), any sequence of equicontinuous functions $X^{c_i} (t)$ on $[0,T]$, $T\in (0, \infty)$, such that $\sup_{c_i} |X^{c_i} (0)| < \infty$, has a converging subsequence with respect to the uniform norm.  In order to apply the theorem, we will show that the sequences
$ (\bar{A}^{c_i}_x, \bar{D}^{c_i}_x,\bar{Q}^{c_i}_x,\bar{\Gamma}_k^{c_i})$ of $4$-tuples
satisfy both properties for $\omega \in G_1$.
\iffalse
\begin{equation}\label{ModConv}
\lim_{\delta \rightarrow \infty} \sup_{X\in S} w_X(\delta) = 0,
\end{equation}
where $w_X(\delta)$ is the modulus of continuity 
\begin{equation*}
w_X(\delta) = \sup_{\substack{0\leq s,t \leq T:\\ |t-s|< \delta }} \left| X(t) - X(s) \right|.
\end{equation*}
\fi

Since $|\bar{\vecQ}^{c_i}(0)|=1$ and the other variables are initially $0$, the second property
clearly holds.  In order to show equicontinuity, we first note that each packet requires one unit of service, and so $\bar{D}^c_x(t)$ is Lipschitz for $x\in\mJ\cup \mR \cup \mK$, $k\in\mK$; by
equation \eqref{eq:gammak}, so is $\bar{\Gamma}^c_k(t)$.
Hence, both variables are equicontinuous. 
%
%By the triangle inequality, for $s,t\in [0,T]$ with $s<t$,
%\begin{equation*}
%\big| \bar{A}^{c_i}_r(t)-\bar{A}_r^{c_i}(s) \big| \leq  2\! \sup_{ 0\leq t\rq{} \leq T} \big| \bar{A}^{c_i}_r(t\rq{}) -a_r t\rq{} \big| + a_r\big| t- s \big|.
%\end{equation*}
%

On the other hand, equicontinuity of $\bar{A}^{c_i}_r (t)$ follows by applying 
the Functional Strong Law of Large Numbers to $\bar{A}^{c_i}_r(t)$ for increasingly fine meshes in $t$. Equicontinuity for the remaining variables $\bar{A}^{c_i}_x (t)$ and $\bar{Q}^{c_i}_x (t)$ follows from straightforward calculations using the equations (\ref{eq:Incr}-\ref{eq:ADk}). Consequently, by the Azel\`a-Ascoli Theorem, every subsequence $\{ (\bar{A}^{c_i}_x, \bar{D}^{c_i}_x,\bar{Q}^{c_i}_x,\bar{\Gamma}_k^{c_i} : x\in \mJ \cup\mK \cup \mR, k\in\mK) \}$ has a further subsequence that converges uniformly on $[0,T]$.

We now show that equations (\ref{eq:Incr}-\ref{eq:ADk}) and (\ref{Fluid:1}-\ref{PF}) hold for 
the limit of any 
converging subsequence of $\{ (\bar{A}^{c_i}_x, \bar{D}^{c_i}_x,\bar{Q}^{c_i}_x,\bar{\Gamma}_k^{c_i} : x\in \mJ \cup\mK \cup \mR, k\in\mK) \}$; we will demonstrate this on $G_1$ for
all equations except (\ref{Fluid:3}-\ref{PF}).
With the exception of (\ref{eq:Dk}-\ref{eq:Ak}), equations (\ref{eq:Incr}-\ref{eq:ADk}) 
clearly hold for any limit since each holds under the prelimit. We must take a little more care for expressions 
(\ref{eq:Dk}-\ref{eq:Ak}),  because some minor error is introduced by linear interpolating terms; the bound
\begin{equation*}
\left| D^{c}_k(s) - \Gamma^{c}_k(D^c_j(s)) \right| \leq \left| D^c_k(\lceil s \rceil) - D^c_k(\lfloor s\rfloor)\right| \leq \sigma_{\max}
\end{equation*}
shows that this error vanishes in the limit. 
On $G_1$, \eqref{Fluid:1} follows directly from the Strong Law of Large Numbers. 

We still need to demonstrate \eqref{Fluid:3} for $\sigma (\vecQ)$ satisfying (\ref{PF}), which we show for $\omega \in G_2$,
where $G_2$ will be chosen as in Lemma \ref{lastlemma}.
Let $(Q_j,D_j : j\in\mJ)$ be the limit of a subsequence of $\{ (\bar{Q}_j^{c_i}, \bar{D}_j^{c_i}: j\in\mJ ) \}$. Suppose that $Q_j(t)>\epsilon$ for some $\epsilon>0$ and a given $j$.  Then, by the continuity of $\vecQ(t)$, there exists $\delta>0$ such that $Q_j(t+h)>\epsilon$ for all $h$ satisfying $|h|\le \delta$. Since convergence is uniform, there also exists $c\rq{}$ such that, for all $c_i\ge c\rq{}$, $\bar{Q}^{c_i}_j(t+h)>\epsilon$ for all such $h$. 

For $u\in\bR_+$, denote by $\lfloor u \rfloor_{c_i}$ be the largest value of $k c_{i}^{-1}$, $k\in\bZ$, that is at most $u$.
%\[
%\lfloor u \rfloor_c := \max_{k\in\bZ_+}\{ k  c^{-1}  :  k c^{-1} \leq u\}.
%\]
%
The scaled departures over $(\lfloor t \rfloor_{c_i}, \lfloor t+h \rfloor_{c_i}]$ can be rewritten as 
\begin{align}
&\bar{D}^{c_i}_j(\lfloor t+h \rfloor_{c_i}) - \bar{D}^{c_i}_j(\lfloor t \rfloor_{c_i})\notag 
= \frac{1}{c_i}\sum_{\tau=\lfloor t \rfloor_{c_i}+1/c_i}^{\lfloor t+h \rfloor_{c_i}} \pi_j(\bar{\vecQ}^{c_i}(\tau))\notag \\
=  &\frac{1}{c_i}\sum_{\tau=\lfloor t \rfloor_{c_i}+1/c_i}^{\lfloor t+h \rfloor_{c_i}} \sigma_j(\bar{\vecQ}^{c_i}(\tau))  + \bar{M}^{c_i}_j(\lfloor t+h \rfloor_{c_i}) - \bar{M}^{c_i}_j(\lfloor t \rfloor_{c_i})\notag \\
= & \int_{\lfloor t \rfloor_{c_i} +1/c_i}^{\lfloor t+h \rfloor_{c_i}} \sigma_j(\bar{\vecQ}^{c_i}(\lfloor \tau \rfloor_{c_i}))   d \tau + \bar{M}^{c_i}_j(\lfloor t+h \rfloor_{c_i}) - \bar{M}^{c_i}_j(\lfloor t \rfloor_{c_i})%+ \frac{1}{c_i} \sum_{\tau=\lfloor t \rfloor_{c_i}+c_i^{-1}}^{\lfloor t + h \rfloor_{c_i}} \left[ \pi_j(\bar{\vecQ}^{c_i}(\tau)) - \sigma_j(\bar{\vecQ}^{c_i}(\tau)) \right]
\label{FinalBound}\,,
\end{align}
%%[MORE DETAIL ON CHANGE OF TIME NEEDED HERE] 
%where $\sigma_j (\cdot)$ and $\pi_j (\cdot)$ are as in (\ref{DServe})  
%
where the summations are understood to be over values in $\{0,c_{i}^{-1}, 2c_{i}^{-1},...\}$.
The square bracketed terms on the last line of \eqref{FinalBound} are bounded.
For $\omega\in G_2$, with $G_2$ as in Lemma \ref{lastlemma}, the last two terms in \eqref{FinalBound} converge uniformly to zero. 
%[CHECK DOOB]
On the other hand, by uniform convergence, the terms $\bar{\vecQ}^{c_i}(\lfloor \tau \rfloor_{c_i})$ on the subsequence converge to $\vecQ_j(\tau)$. Also, $\sigma_j(\vecQ)$ is bounded and continuous on $Q_j>0$. Thus, applying the Bounded Convergence Theorem, 
it follows from (\ref{FinalBound}) that
\begin{equation*}
{D}_j(t+h) - {D}_j(t) = \int_{t}^{t+h} \sigma_j({\vecQ}(\tau))   d \tau\,,
\end{equation*}
where $\sigma_j$ is as in (\ref{PF}).
Dividing by $h$ and taking the limit $h\searrow 0$ implies $D\rq{}_j(t) = \sigma_j(\vecQ(t))$, and hence that \eqref{Fluid:3} also holds.
\Halmos\endproof

%\begin{lemma} 
%\label{lastlemma}
%Almost surely, for $j\in\mJ$ and for all $T>0$
%\[
%\sup_{t\leq T} \Big|\bar{M}_j^c(t)  \Big|\xrightarrow[c\rightarrow \infty]{} 0.
%\]
%\end{lemma}
%
%

\smallskip

\proof{Proof of Lemma \ref{lastlemma}}
%%For each $c$ , we define 
%%%% CHECK t FLOOR
%%\[
%%M^c(t) = \sum_{\tau=1}^{\lfloor t\rfloor_c} \left[ \pi_j ( \bar{\vecQ}^c(\tau)) - \sigma_j( \bar{\vecQ}^c(\tau)) \right],
%%\]
%%where $t \in \{ 0, c^{-1}, 2c^{-1},...\}$ and where, for this proof, summations and \NW{maximizations} are understood to be over values of \NW{$s$, $t$ and $T$ in $\{0,c^{-1}, 2c^{-1},...\}$}.  
%
For $\epsilon > 0$, let
$\tau_\epsilon$ denote the stopping time for the event $\{ |M_j^c(t)| \ge c\epsilon  \}$.  
Since, for each $c$, $M_j^c$ is a martingale with increments uniformly bounded by 
$\sigma_{\text{max}}$, it follows by applying 
the Azuma-Hoeffding Inequality \cite[E14.2]{Wi91} to $\{ M_j^{c}(\min(\tau_\epsilon, s)) \}_{s\leq cT}$ that
\begin{align*}
  \bP \left(  \max_{t: t \leq c T} \left| M_j^{c}(t)\right| \ge c\epsilon  \right) 
= & \bP\Big( \big|M_j^c\big(\min(\tau_\epsilon, c \, T)\big)\big| \ge c\epsilon \Big) \leq  2 \exp\left\{ - \frac{c^2\epsilon^2}{2 cT \sigma_{\max}^2 } \right\}.
\end{align*}
After summing over $c\in\bN$, the right hand side of the above inequality is finite and hence,
by the Borel-Cantelli Lemma, 
\[
\limsup_{c\rightarrow \infty}  \max_{t: t\leq c T} \frac{1}{c} \left| M_j^{c}(t) \right|  \leq \epsilon
\]
holds almost surely.  Since $\epsilon >0$ is arbitrary, this implies the claim in the lemma. \Halmos
\endproof

\iffalse
\NW{We briefly comment that if we allowed for schedules to be random, as discussed in Remark \ref{ScheduleRemark}, then the convergence of   \eqref{FinalBound} would still hold under the Martingale Convergence Theorem. And thus fluid model under random schedules (with expected values $\mS$) is given by the same fluid model equations. Since our fluid analysis does not require the set of schedules to be integral, positive recurrence holds in this case also.}
\fi

\section{Positive Recurrence}
\label{Appendix F}

In this subsection, we demonstrate Proposition \ref{propFMQN}, which states that 
positive recurrence for a FIFO proportional switched network follows
from the fluid model stability of the corresponding  proportional switch fluid model.  The
main tools are Proposition \ref{FluidLimit} of the previous subsection and the Multiplicative Foster\rq{}s Criterion.

\proof{Proof of Proposition \ref{propFMQN}}
%We follow the approach due to \cite{Da95} and surveyed in \cite{Br08}.
%Theorem 4.16 of \cite{Br08} states that, under appropriate regularity conditions that include the arrival process at each queue being a renewal process, a queueing 
%network will be positive recurrent if its corresponding fluid model is
%stable.  We cannot directly apply Theorem 4.16 to our setting since our
%discrete time queueing network has input that is i.i.d. at each time but is not a renewal process.  However, Theorem \ref{mainthrm} will follow immediately from Theorem \ref{FluidStable} (which shows stability of the fluid model) and a minor variant of Theorem 4.16.

%Given our fluid stability result, Theorem \ref{FluidStable}, positive recurrence follows from standard stability surveyed by \cite{Br08} in Theorem 4.16 and for FIFO queueing networks in Theorem 5.9. The only minor differences being the discrete time setting, our model is a countable state-space Markov chain (thus recurrence concepts simplify, see \cite{Br08} page 280), arrivals are compound renewal processes rather than renewal process.

%The key aspect of the proof of fluid stability implying positive recurrence is the reduction to Multiplicative Foster\rq{}s Criterion \cite[Proposition 4.6]{Br08}, which we now describe. 
%%The state space $\mathcal X$ of the Markov process $(\vecX(t), t\geq0)$ has been defined after \eqref{eq:Markov.Chain.X}.
%%Let define the modulus of the state of the Markov process in the following way
%%$$||\vecX (t)|| = |\vecQ(t)| \ .$$

One can check that, for each $t\geq 0$, the sequence of queue sizes 
$\{ |\bar{\vecQ}^{c}(t)| \}_c$ of the
scaled proportional switched networks in Proposition \ref{FluidLimit}
is uniformly integrable. This follows quickly from the inequality
\begin{equation}
\label{equnifint}
|\bar{\vecQ}^{c}(t)| =\sum_{j\in\mJ} \bar{Q}_j^{c}(t)\leq \sum_{j\in\mJ} \frac{Q^{c}_j(0)}{c} + \sum_{r\in\mR} \frac{A^{c}_r(t)}{c}
\end{equation}
since $A^{c}_r(t)$ is a sum of i.i.d. random variables with finite mean (see, e.g., \citet[ Lemma 4.13, (4.81)]{Br08}). %If more details needed apply union bound to definition of uniform ingrability

On the other hand, by Proposition \ref{FluidLimit}, on a set of probability one, every subsequence of
$\bar{\vecQ}^{c}(t)$ has a further subsequence that converges uniformly on compact time
 intervals to a fluid model solution $\vecQ(t)$ of (\ref{eq:Incr}--\ref{eq:ADk}) and (\ref{Fluid:1}--\ref{PF}), with $|\vecQ(0)|=1$. By Theorem \ref{FluidStable}, this fluid model is stable; hence all fluid model solutions with $|\vecQ(0)|=1$  satisfy $|\vecQ(t)| = 0$ for $t\ge \gamma$, with $\gamma$ not depending on the particular
fluid model solution. 

It follows from this that, on a set of probability one, every subsequence of 
$|\bar{\vecQ}^{c}(\gamma)|$ has a further subsequence converging to zero; consequently,  
$|\bar{\vecQ}^{c}(\gamma)|$ also converges to zero on a set of probability one. By the above uniform integrability and almost sure convergence of $|\bar{\vecQ}^{c}(\gamma)|$, it follows that
\begin{equation*}
\lim_{c\rightarrow\infty} \bE |\bar{\vecQ}^{c}(\gamma)|  =0\,,
\end{equation*}
and so, if $c\ge \tilde{c}$ for appropriate $\tilde{c}$, then
$\bE |\bar{\vecQ}^{c}(\gamma)|   \le 1/2$.  

The last inequality is equivalent to 
$\bE |\vecQ^{c}(\gamma c)|   \le c/2$; together with (\ref{equnifint}), it implies that, for large enough $\kappa$ and all $c$, 
\begin{equation}
\label{eqMFC}
\bE |\vecQ^{c}(\gamma (\max(c,\kappa)))|   \le \max(c,\kappa)/2 \,.
\end{equation}
%
%By definition $c=|\vecQ^{c}(0)|$, and $\bar{\vecQ}^{c}(t)= \vecQ^{c}(t)/c$ given a process started from $\vecQ^{c}(0)$ with $|\vecQ^{c}(0)|=c$. The above inequality can be rearranged and states for a (discrete-time) FIFO switched network under the Proportional Scheduler, $\tilde{\vecQ}$, the following holds 
%\begin{equation*}
%\bE\Big[  \big|{\tilde{\vecQ}}(T^*|\tilde{\vecQ}(0)|)\big| - \big|\tilde{\vecQ}(0)\big| \Big| \tilde{\vecQ}(0) \Big]  < -\epsilon |\tilde{\vecQ}(0)|
%\end{equation*}
%when $|\tilde{\vecQ}(0)|>\tilde{c}$. 
%
The inequality (\ref{eqMFC}) satisfies the main premise of the Multiplicative Foster\rq{}s Criterion \cite[Proposition 4.6, (4.28)]{Br08}, and therefore implies the positive recurrence of the corresponding
proportional switched network, which completes the proof of the proposition.  
(The criterion is stated in \cite{Br08} for continuous time
Markov processes; however, both the criterion and its proof carry over to discrete time.  Alternatively, the
discrete time Markov process can be embedded in continous time.  The petite set assumption in the criterion is automatically satisfied in our framework since the empty state will be hit with uniformly high
probability from sets with a bounded number of packets.)
\iffalse
Since the state space is discrete, we have that for any $K>0$ the set
$$F_K = \{ x \in \mathcal X: ||x|| \leq K\}$$
is finite.
From \eqref{eq:Q.hitting.time.of.zero} and \eqref{eq:fluid.scaled.processes} we have that, with $T=c_5$
$$
\limsup_{\vecQ(0) \to\infty} \bE_{\vecQ(0)}[||\bar X(T)||] 
= \limsup_{\vecQ(0)\to\infty} \frac{\bE_{\vecQ(0)}[||\bar X(|\vecQ(0)| \, T)||] }{|\vecQ(0)|} = 0 \ .
$$
The property follows from a derivate result of the Foster's criterion, see Corollary 9.8, page 259 in \cite{Ro03}.
\fi
\Halmos\endproof

\end{document}